# Near-Optimal Stochastic Approximation for Online Principal Component Estimation

Chris Junchi Li *    Mengdi Wang *    Han Liu *    Tong Zhang †

October 5, 2017

**Abstract**

Principal component analysis (PCA) has been a prominent tool for high-dimensional data analysis. Online algorithms that estimate the principal component by processing streaming data are of tremendous practical and theoretical interests. Despite its rich applications, theoretical convergence analysis remains largely open. In this paper, we cast online PCA into a stochastic nonconvex optimization problem, and we analyze the online PCA algorithm as a stochastic approximation iteration. The stochastic approximation iteration processes data points incrementally and maintains a running estimate of the principal component. We prove for the first time a nearly optimal finite-sample error bound for the online PCA algorithm. Under the subgaussian assumption, we show that the finite-sample error bound closely matches the minimax information lower bound.

**Keywords:** Principal component analysis, Stochastic approximation, Nonconvex optimization, Stochastic gradient method, High-dimensional data, Online algorithm, Finite-sample analysis.

## 1 Introduction

Principal component analysis (PCA) (Pearson, 1901; Hotelling, 1933) is one of the most popular dimension reduction methods for high-dimensional data analysis. It has wide applications in bioinformatics, healthcare, imaging, computer vision, artificial intelligence, social science, finance and economy. Let $\boldsymbol{X}$ be a random vector in $\mathbb{R}^d$ with mean zero and unknown covariance matrix

$$\boldsymbol{\Sigma} = \mathbb{E}\left[\boldsymbol{X}\boldsymbol{X}^\top\right] \in \mathbb{R}^{d\times d},$$

where the eigenvalues of $\boldsymbol{\Sigma}$ are $\lambda_1 > \lambda_2 \geq \cdots \geq \lambda_d \geq 0$. Principal component analysis aims to find the principal eigenvector of $\boldsymbol{\Sigma}$ that corresponds to the largest eigenvalue $\lambda_1$, based on independent and identically distributed sample realizations $\boldsymbol{X}^{(1)},\ldots,\boldsymbol{X}^{(n)}$. This can be casted into a *nonconvex*

---

*Department of Operations Research and Financial Engineering, Princeton University, Princeton, NJ 08544, USA; e-mail: {junchil, mengdiw, hanliu}@princeton.edu

†Tencent AI Lab, Shennan Ave, Nanshan District, Shenzhen, Guangdong Province 518057, China; email: tongzhang@tongzhang-ml.org

*stochastic optimization problem*, given by

$$\begin{aligned}
\text{maximize} \quad & \mathbf{u}^\top \mathbb{E}\left[\boldsymbol{X}\boldsymbol{X}^\top\right]\mathbf{u}, \\
\text{subject to} \quad & \|\mathbf{u}\| = 1, \\
& \mathbf{u} \in \mathbb{R}^d,
\end{aligned} \quad (1.1)$$

where $\|\cdot\|$ denotes the Euclidean norm. We assume throughout this paper that *the covariance matrix $\boldsymbol{\Sigma}$ has a unique principal component*, which we denote by $\mathbf{u}^*$. Accordingly, the principal component $\mathbf{u}^*$ is the unique solution to problem (1.1).

The classical PCA method estimates $\mathbf{u}^*$ using the principal component of the empirical covariance matrix, i.e.,

$$\widehat{\mathbf{u}}^{(n)} = \operatorname*{argmax}_{\|\mathbf{u}\|=1} \mathbf{u}^\top \widehat{\boldsymbol{\Sigma}}^{(n)} \mathbf{u},$$

where $\widehat{\boldsymbol{\Sigma}}^{(n)}$ is the empirical covariance matrix based on $n$ samples

$$\widehat{\boldsymbol{\Sigma}}^{(n)} = \frac{1}{n} \sum_{i=1}^{n} \boldsymbol{X}^{(i)} \left(\boldsymbol{X}^{(i)}\right)^\top.$$

This can be viewed as a *sample average approximation* method for problem (1.1). According to both literatures of PCA and sample average approximation, it is well-known that this method produces a *non-improvable* solution $\widehat{\mathbf{u}}^{(n)}$, in the sense that the estimation error $\widehat{\mathbf{u}}^{(n)} - \mathbf{u}^*$ achieves the information lower bound. More precisely, we are interested in the angle between the two unit vectors $\widehat{\mathbf{u}}^{(n)}$ and $\mathbf{u}^*$, given by

$$\angle(\widehat{\mathbf{u}}^{(n)}, \mathbf{u}^*) = \arccos\left(\widehat{\mathbf{u}}^{(n)\,\top} \mathbf{u}^*\right), \quad (1.2)$$

which takes value in $[0, \pi]$. It has been shown that the angle between any $n$-sample estimator $\widetilde{\mathbf{u}}^{(n)}$ and the true principal component $\mathbf{u}^*$ satisfies

$$\inf_{\widetilde{\mathbf{u}}^{(n)}} \sup_{\boldsymbol{X} \in \mathcal{M}(\sigma_*^2, d)} \mathbb{E}\left[\sin^2 \angle(\widetilde{\mathbf{u}}^{(n)}, \mathbf{u}^*)\right] \geq c \cdot \sigma_*^2 \cdot \frac{d-1}{n}, \quad (1.3)$$

where $c$ is some positive constant (Theorem 3.1 of Vu & Lei (2013)). In Eq. (2.18), the infimum of $\widetilde{\mathbf{u}}^{(n)}$ is taken over all $n$-sample estimators, and $\mathcal{M}(\sigma_*^2, d)$ is the collection of all $d$-dimensional subgaussian distributions with mean zero and eigenvalues satisfying $\lambda_1 \lambda_2 / (\lambda_1 - \lambda_2)^2 \leq \sigma_*^2$. Classical PCA method has time complexity $\mathcal{O}(nd^2)$ and space complexity $\mathcal{O}(d^2)$. When the raw data are high-dimensional, storing and computing a large empirical covariance matrix can be expensive.

In this paper, we focus on *online methods* for principal component analysis. These methods update the iterates incrementally by processing data points one by one. The practical goal is to be able to learn the principal eigenvector "on the fly", without explicitly computing and storing the empirical covariance matrix $\widehat{\boldsymbol{\Sigma}}$. In particular, we focus on an iteration for online PCA that was first proposed by Oja (1982), which is given by

$$\mathbf{u}^{(n)} = \Pi \left\{ \mathbf{u}^{(n-1)} + \beta \boldsymbol{X}^{(n)} (\boldsymbol{X}^{(n)})^\top \mathbf{u}^{(n-1)} \right\}, \quad (1.4)$$

where $\beta$ is some positive stepsize, and $\Pi$ denotes the Euclidean projection operator onto the unit sphere $\mathcal{S}^{d-1} = \{\mathbf{u} \in \mathbb{R}^d \mid \|\mathbf{u}\| = 1\}$, i.e., $\Pi\mathbf{u} = \|\mathbf{u}\|^{-1}\mathbf{u}$ for all $\mathbf{u} \neq 0$. The iteration (1.4) only requires vector product operations. It has time complexity $\mathcal{O}(d)$ per iteration, and has space complexity $\mathcal{O}(d)$. Iteration (1.4) is very easy to implement in practice and has been used as a heuristic method for fast principal component analysis.

In contrast to classical PCA which is a sample average approximation method, the online PCA iteration (1.4) is essentially a *stochastic approximation* iteration for the optimization problem (1.1). Surprisingly, theoretical analysis of its convergence remains largely open (except in special cases). The theoretical challenge is due to our attempt to maximize a convex function over a unit sphere, on which there are infinitely many unstable stationary solutions. In this paper, we aim to take an important step in analyzing the convergence of online PCA methods. We are interested in finite-sample analysis of the online PCA iteration, and we aim to match the convergence rate with the information lower bound. An intriguing open question is: *Can online principal component analysis be optimal?* An answer "yes" would imply that principal component analysis is as simple as estimating the mean of a distribution, for which an iterative online process is sufficient to provide non-improvable estimates. Our results provide a partial answer to this question: *almost yes*. Our analysis involves the weak convergence theory for Markov processes (Ethier & Kurtz, 2005). Such analysis has a potential for a broader class of stochastic algorithms for nonconvex optimization.

## 1.1 Related Literatures

Although the online PCA iteration (1.4) was proposed over thirty years ago (Oja, 1982; Oja & Karhunen, 1985), its convergence rate analysis remains somewhat limited. It was not until recently, due to the need to handle massive amounts of data, did online PCA gain attention. Two recent works Balsubramani et al. (2013) and Shamir (2016a) study the convergence of online PCA from different perspectives, and obtain some useful rate results. We provide a detailed comparison between their results and ours in §2.3.

In the mathematical programming and statistics communities, the computational and statistical aspects of PCA are often studied separately. From the statistical perspective, recent developments have focused on estimating principal components for very high-dimensional data. When the data dimension is much larger than the sample size, i.e., $d \gg n$, classical method using decomposition of the empirical covariance matrix produces inconsistent estimates (Johnstone & Lu, 2009; Nadler, 2008). Sparsity-based methods have been studied, such as the truncated power method studied by Yuan & Zhang (2013) and Wang et al. (2014). Other sparsity regularization methods for high dimensional PCA has been studied in Johnstone & Lu (2009); Vu & Lei (2013, 2012); Zou (2006); d'Aspremont et al. (2008); Amini & Wainwright (2009); Ma (2013); Cai et al. (2013), etc. Note that in this paper we do not consider the high-dimensional regime and sparsity regularization.

From the computational perspective, power iterations or the Lanczos method are well studied. These iterative methods require performing multiple products between vectors and empirical covariance matrices. Such operation usually involves multiple passes over the data, whose complexity may scale with the eigengap and dimensions (Kuczynski & Wozniakowski, 1992; Musco & Musco, 2015). Recently, randomized algorithms have been developed to reduce the computation complexity (Shamir, 2015, 2016b; Garber & Hazan, 2015). A critical trend today is to combine the computational and statistical aspects and to develop algorithmic estimator that admits fast computation as well as good estimation properties. Related literatures include Arora et al. (2012, 2013); Mitliagkas et al. (2013); Hardt & Price (2014); De Sa et al. (2015).

The idea of using stochastic approximation for PCA can be traced back to Oja & Karhunen (1985). Stochastic approximation (SA) was first studied for the root finding problem and later extended to stochastic optimization and stochastic variational inequalities; see e.g., the textbooks by Kushner & Yin (2003), by Benveniste et al. (2012), by Borkar (2008), by Bertsekas & Tsitsiklis (1989). The idea of processing one sample at a time is also related to the class of incremental methods, which are developed for minimizing the sum of a large number of component functions. These methods update incrementally by making use of one component at a time, through a gradient-type or proximal-type iteration [see for example, Nedić et al. (2001); Nedić & Bertsekas (2001); Bertsekas (2011); Nedić (2011); Wang & Bertsekas (2014b,a); Wang et al. (2017)]. However, existing convergence rate analysis of incremental algorithms do not apply to the optimization over the nonconvex sphere constraint.

In the context of machine learning and signal processing, stochastic approximation (or more commonly referred to as stochastic gradient descent) has been extensively studied for stochastic convex optimization. In machine learning applications, the objective is usually the expectation of a convex loss function parameterized by a random variable. It has been shown that after $n$ samples/iterations, the average of the iterates has $\mathcal{O}(1/n)$ optimization error for strongly convex objective, and $\mathcal{O}(1/\sqrt{n})$ error for general convex objective (see Rakhlin et al. (2012); Shamir & Zhang (2013)). For nonsmooth problems with noisy gradients, there are $\Theta(1/n)$ and $\Theta(1/\sqrt{n})$ minimax information-theoretic lower bounds for convex and strongly convex problems, respectively (see e.g., Agarwal et al. (2012), and see the book by Nemirovsky & Yudin (1983) for a comprehensive study on optimization complexities). The matching between the lower bounds and the upper bounds suggests that stochastic approximation (or equivalently, stochastic gradient) is *optimal* for convex optimization under the stochastic first-order oracle. In contrast, to the best knowledge of the authors, there has been no such work for stochastic nonconvex optimization problem (1.1) with sphere constraint.

## 1.2 Our Contributions

The contributions of this work is summarized as follows.

- We provide the first convergence rate result for online PCA that *nearly matches the information minimax lower bound* (2.18) under the subgaussian distributional assumption. We show that, when the initial iterate $\mathbf{u}^{(0)}$ is randomly chosen according to a uniform distribution and the stepsize $\beta$ is chosen in accordance with the sample size $n$, there is a high-probability event $\mathcal{A}_*$ with $\mathbb{P}(\mathcal{A}_*) \geq 1 - \delta$ such that

$$\mathbb{E}\left[\sin^2 \angle(\mathbf{u}^{(n)}, \mathbf{u}^*) \mid \mathcal{A}_*\right] \leq C(d, n, \delta) \cdot \frac{\lambda_1 \lambda_2}{(\lambda_1 - \lambda_2)^2} \cdot \frac{(d-1)\log n}{n}, \tag{1.5}$$

where $\delta \in [0, 1)$ and $C(d, n, \delta)$ is some factor that can be approximately treated as a constant; see details in Eq. (2.16). We also show that, when both the data dimension $d$ and data size $n$ scale up with $d/n^{1-\varepsilon} \to 0$ for some constant $\varepsilon \in (0,1)$, the factor $C(d,n,\delta)$ approaches to some absolute constant. Moreover, we show that

$$\mathbb{E}\left[\sin^2 \angle(\mathbf{u}^{(n)}, \mathbf{u}^*)\,;\, \mathcal{A}_*\right] \leq C'(d,n,\delta) \cdot \frac{\lambda_1}{\lambda_1 - \lambda_2} \sum_{k=2}^{d} \frac{\lambda_k}{\lambda_1 - \lambda_k} \cdot \frac{\log n}{n}, \tag{1.6}$$

where the factor $C'(d,n,\delta)$ converges to some absolute constant as long as $d^{-1}\sum_{k=1}^{d}(\lambda_k/\lambda_1)$ is bounded away from 0. To the authors' best knowledge, this is the tightest convergence result known for the online PCA iteration under the *near-optimal scaling condition*, i.e., $d/n^{1-\varepsilon} \to 0$. A detailed comparison between our result and concurrent ones is given in §2.3.

- Our convergence rate results are *nearly global*. In particular, convergence rate results hold as long as the initial iterate satisfies

$$|\sin \angle(\mathbf{u}^{(0)}, \mathbf{u}^*)| \leq 1 - \frac{c}{d},$$

for some constant $c > 0$. Here the error tolerance scales inverse-linearly as the dimension $d$ increases. This is critical because, when $d$ is large, a uniformly distributed initial iterate is *nearly perpendicular* to the principal component with high probability. Our initial condition allows one to randomly sample $\mathbf{u}^{(0)}$ according to a uniform distribution over the sphere, while preserving the near-optimal convergence rate. In contrast, most existing results on PCA, for instance Yuan & Zhang (2013), require that the initial condition be $|\sin \angle(\mathbf{u}^{(0)}, \mathbf{u}^*)| \leq 1 - c_0$ for some $c_0 \in (0,1)$. Such initial condition becomes increasingly stringent as the dimension $d$ increases (Ball, 1997). Our choice of initial iterate does not require any prior knowledge about the principal component. Therefore, our convergence results are nearly global in the sense that a randomly selected initial point achieves near-optimal convergence rate with high probability.

**Organization.** §2 states the assumptions on subgaussian PCA and give the main results on convergence rate and finite-sample error bound of the online PCA iteration. §3 analyzes the convergence of the stochastic iteration and gives proofs of the main results. §4 summarizes the results of this paper and presents directions for future works. Some technical analyses are deferred to the Supplementary Material.

**Notations.** For a sequence of $\{x_k\}$ and positive $\{y_k\}$, we write $x_k = \mathcal{O}(y_k)$ if there exists constant $M < \infty$ such that $|x_k| \leq My_k$. Also let $\lfloor x \rfloor$ denotes the largest integer $\leq x$, and $\lceil x \rceil$ the smallest integer $\geq x$. Let $x \wedge y = \min(x,y)$ and $x \vee y = \max(x,y)$. Let $\|\mathbf{A}\|_2$ be the spectral norm of square matrix $\mathbf{A}$. Lastly, let $\|X\|_\infty$ be the $L^\infty$-norm of random variable $X$, and also $\mathbb{E}[X; \mathcal{A}]$ denotes the expectation of $X$ over event $\mathcal{A}$, i.e. we have $\mathbb{E}[X; \mathcal{A}] = \mathbb{E}[X\mathbf{1}_\mathcal{A}]$.

## 2 Main Results

In this section, we present the main convergence results for the online PCA algorithm. The algorithm maintains a running estimate $\mathbf{u}^{(n)}$ of the true principal component $\mathbf{u}^*$, and updates it while interacting with a streaming data source. Due to the nature of PCA, we are interested in the convergence of the angle process $\{\angle(\mathbf{u}^{(n)}, \mathbf{u}^*)\}$. A detailed description of the algorithm is given by Algorithm 1.

**Algorithm 1** Online PCA Algorithm
---
Initialize $\mathbf{u}^{(0)}$ and choose the stepsize $\beta$ under some condition (to be specified)
**for** $n = 1, 2, \ldots, N$ **do**
  Draw one sample $\boldsymbol{X}^{(n)}$ from the (streaming) data source
  Update the iterate $\mathbf{u}^{(n)}$ by
$$\mathbf{u}^{(n)} = \Pi \left\{ \mathbf{u}^{(n-1)} + \beta \boldsymbol{X}^{(n)} (\boldsymbol{X}^{(n)})^\top \mathbf{u}^{(n-1)} \right\}.$$
**end for**
---

## 2.1 Distributional Assumptions

In this paper we focus on the setting where $\boldsymbol{X}$ has a subgaussian distribution. We first introduce the subgaussian norm of a random variable, following Vershynin (2010).

**Definition 1** (Subgaussian norm). For a real-valued random variable $Y$, we define its subgaussian norm as

$$\|Y\|_{\psi_2} \equiv \sup_{p \geq 1} p^{-1/2} \left(\mathbb{E}|X|^p\right)^{1/p}. \tag{2.1}$$

When a zero-mean random variable $Y$ has a finite subgaussian norm, its distribution is known to have a *subgaussian tail*. On the other hand, if $Y \sim N(0, \sigma^2)$, we can verify that $\|Y\|_{\psi_2} \leq \sigma\sqrt{2/\pi}$.

**Definition 2** (Vector-subgaussian). We say that a random vector $\boldsymbol{Z} \in \mathbb{R}^d$ follows a subgaussian distribution if $\sup_{\mathbf{v} \in \mathcal{S}^{d-1}} \|\mathbf{v}^\top \boldsymbol{Z}\|_{\psi_2} < \infty$.

We refer the readers to §B and van der Vaart & Wellner (1996); Vershynin (2010) for more discussions on subgaussian distributions. Let $\boldsymbol{\Sigma}^{1/2}$ be the unique symmetric positive semidefinite matrix that satisfies $\boldsymbol{\Sigma}^{1/2} \boldsymbol{\Sigma}^{1/2} = \boldsymbol{\Sigma}$. We state our first assumption as follows.

**Assumption 1** (Subgaussian Distribution). Let $\boldsymbol{Z} \in \mathbb{R}^d$ be a random vector satisfying

$$\mathbb{E}[\boldsymbol{Z}] = \mathbf{0}, \quad \mathbb{E}\left[\boldsymbol{Z}\boldsymbol{Z}^\top\right] = \mathbf{I}_d, \quad \sup_{\mathbf{v} \in \mathcal{S}^{d-1}} \|\mathbf{v}^\top \boldsymbol{Z}\|_{\psi_2} \leq 1. \tag{2.2}$$

Let $\boldsymbol{X}^{(1)}, \boldsymbol{X}^{(2)}, \ldots \in \mathbb{R}^d$ be independent and identically distributed realizations of $\boldsymbol{X} \equiv \boldsymbol{\Sigma}^{1/2} \boldsymbol{Z}$.

Assumption 1 essentially requires that streaming data $\boldsymbol{X}^{(1)}, \boldsymbol{X}^{(2)}, \ldots$ follow a subgaussian distribution when projected onto every unit vector in $\mathbb{R}^d$. It follows immediately from Assumption 1 that $\mathbb{E}[\boldsymbol{X}] = 0$ and $\mathbb{E}\left[\boldsymbol{X}\boldsymbol{X}^\top\right] = \boldsymbol{\Sigma}$. It is also straightforward to verify that the normal distribution $\boldsymbol{X}$ with mean $\mathbf{0}$ and covariance $\boldsymbol{\Sigma}$ satisfies Assumption 1.

**Assumption 2** (Positive Eigenvalue Gap). Let the eigenvalues of $\boldsymbol{\Sigma}$ be $\lambda_1 > \lambda_2 \geq \ldots \geq \lambda_d \geq 0$.

The positive value $\lambda_1 - \lambda_2$ is often referred to as the *eigengap*. The nonzero eigengap implies that the principal component $\mathbf{u}^*$ is uniquely identified. We will show later that the eigengap plays a key role in the convergence analysis of the online PCA iteration.

## 2.2 Convergence Results and Finite Sample Analysis

For easiness of presentation, we introduce a rescaling of the iteration/sample index $n$. We define the rescaled time $N^*_{\beta,s}$ as

$$N^*_{\beta,s} = \left\lceil \frac{s \log(\lambda_1^{-2}(\lambda_1 - \lambda_2)\beta^{-1})}{-\log(1 - \beta(\lambda_1 - \lambda_2))} \right\rceil, \tag{2.3}$$

where $s > 0$ is a tuning parameter. Here the fixed time $N^*_{\beta,s}$ increases to infinity as the stepsize $\beta$ decreases to 0. In later analysis, we will choose the stepsize $\beta$ to be inversely proportional to the budget sample size; see Eq. (2.11). For any $c^* > 0$, we also define the rescaled time $N^o_\beta(c^*)$ as

$$N^o_\beta(c^*) = \left\lceil \frac{\log(4c^*d)}{-\log(1 - \beta(\lambda_1 - \lambda_2))} \right\rceil. \tag{2.4}$$

We state our first result as follows. It provides a foundation to our finite-sample convergence rate analysis.

**Theorem 1** (Convergence result with deterministic initialization). Let Assumptions 1 and 2 hold, and let there be some constant $c^* > 1$ such that $\tan^2 \angle(\mathbf{u}^{(0)}, \mathbf{u}^*) \leq c^*d$. Then for any $\varepsilon \in (0, 1/8)$ and $\beta > 0$ satisfying

$$d[\lambda_1^2(\lambda_1 - \lambda_2)^{-1}\beta]^{1-2\varepsilon} \leq b_1/c^*, \tag{2.5}$$

and any $t > 1$, there exists an event $\mathcal{H}_*$ with

$$\begin{aligned}\mathbb{P}(\mathcal{H}_*) \geq &\ 1 - 2(d+2)N^o_\beta(c^*) \exp\left(-C_0[\lambda_1^2(\lambda_1 - \lambda_2)^{-1}\beta]^{-2\varepsilon}\right) \\ &- 4dN^*_{\beta,t} \exp\left(-C_1[\lambda_1^2(\lambda_1 - \lambda_2)^{-1}\beta]^{-2\varepsilon}\right),\end{aligned} \tag{2.6}$$

such that the iterates generated by Algorithm 1 satisfy the following for $n \in [N^*_{\beta,1} + N^o_\beta(c^*), N^*_{\beta,t}]$:

$$\begin{aligned}\mathbb{E}\left[\tan^2 \angle(\mathbf{u}^{(n)}, \mathbf{u}^*); \mathcal{H}_*\right] \leq &\ (1 - \beta(\lambda_1 - \lambda_2))^{2(n - N^o_\beta(c^*))} \\ &+ C_2 \sum_{k=2}^d \frac{\lambda_1 \lambda_k + \lambda_1^2 \left[\lambda_1^2(\lambda_1 - \lambda_2)^{-1}\beta\right]^{0.5 - 4\varepsilon}}{\lambda_1 - \lambda_k} \cdot \beta.\end{aligned} \tag{2.7}$$

Here $b_1 \in (0, \ln^2 2/16)$, $C_0$, $C_1$ and $C_2$ are positive absolute constants.

Theorem 1 gives a convergence rate estimate for the online PCA algorithm when the initial condition is nearly global, i.e., $\tan^2 \angle(\mathbf{u}^{(0)}, \mathbf{u}^*) \leq c^*d$ for some $c^* > 0$. It characterizes the time/sample size $n$ needed for the angles of iterates to become sufficiently small. Note that the error bound (2.7) is satisfied on a set $\mathcal{H}_*$ with probability close to 1 even if $d \to \infty, \beta \to 0$, as long as $d\beta^{1-2\varepsilon} \to 0$ (assuming all other parameters fixed). This means that the convergence rate result is useful in the regime of high-dimensional data analysis. The proof of Theorem 1 is provided in §3.1.

Next consider the case where the initial solution $\mathbf{u}^{(0)}$ is sampled from the unit sphere according to a uniform distribution. In this case, we would expect that $\tan^2 \angle(\mathbf{u}^{(0)}, \mathbf{u}^*) = \Theta(d)$ with high probability. This implies the initial condition required by Theorem 1 is indeed satisfied. Our second result is given below.

**Theorem 2** (Convergence result with uniformly randomized initialization)**.** Let Assumptions 1 and 2 hold, let $\mathbf{u}^{(0)}$ be uniformly sampled from the unit sphere. Then for any $\varepsilon \in (0, 1/8)$, $\beta > 0$, $\delta > 0$ satisfying

$$d[\lambda_1^2(\lambda_1 - \lambda_2)^{-1}\beta]^{1-2\varepsilon} \leq b_2\delta^2, \tag{2.8}$$

and

$$4dN^*_{\beta,2} \exp\left(-C_3[\lambda_1^2(\lambda_1 - \lambda_2)^{-1}\beta]^{-2\varepsilon}\right) \leq \delta, \tag{2.9}$$

there exists an event $\mathcal{A}_*$ with $\mathbb{P}(\mathcal{A}_*) \geq 1 - 2\delta$ such that the iterates generated by Algorithm 1 satisfy

$$\mathbb{E}\left[\tan^2 \angle(\mathbf{u}^{(n)}, \mathbf{u}^*) \, ; \, \mathcal{A}_*\right] \leq C_4 \delta^{-4} d^2 \left(1 - \beta(\lambda_1 - \lambda_2)\right)^{2n}$$
$$+ C_4 \sum_{k=2}^{d} \frac{\lambda_1 \lambda_k + \lambda_1^2 \left[\lambda_1^2(\lambda_1 - \lambda_2)^{-1}\beta\right]^{0.5-4\varepsilon}}{\lambda_1 - \lambda_k} \cdot \beta. \tag{2.10}$$

for $n \in [N^*_{\beta,2}, N^*_{\beta,3}]$. Here $b_2$, $C_3$ and $C_4$ are positive absolute constants.

Now let us consider the choice of stepsize $\beta$. Suppose that the eigengap $\lambda_1 - \lambda_2$ is known in advance. Also suppose that we are given a budget of sample size $N$. Our goal is to choose an appropriate constant stepsize $\beta$ and minimize the finite-sample error bound (2.10). We will pick the stepsize to be

$$\bar{\beta}(N) = \frac{2 \log N}{(\lambda_1 - \lambda_2)N} \tag{2.11}$$

which is *asymptotically* the minimizer of the right hand of Eq. (2.10). Then we obtain the main result of this paper: a finite-sample error bound for online PCA under subgassian assumption. It is stated as follows.

**Theorem 3** (Main Result: Finite-sample error bound with uniformly randomized initialization)**.** Let Assumptions 1 and 2 hold, let $\mathbf{u}^{(0)}$ be uniformly sampled from the unit sphere, and let $\beta = \bar{\beta}(N)$ be given by Eq. (2.11). Then for any $\varepsilon_0 \in (0, 1/8)$, $N \geq 1$, $\delta > 0$ satisfying

$$d\left[\lambda_1^2(\lambda_1 - \lambda_2)^{-1}\bar{\beta}(N)\right]^{1-2\varepsilon_0} = d\left[\frac{\lambda_1^2}{(\lambda_1 - \lambda_2)^2} \cdot \frac{2 \log N}{N}\right]^{1-2\varepsilon_0} \leq b_3 \delta^2, \tag{2.12}$$

and

$$4dN^*_{\bar{\beta}(N),2} \exp\left(-C_6[\lambda_1^2(\lambda_1 - \lambda_2)^{-1}\bar{\beta}(N)]^{-2\varepsilon_0}\right) \leq \delta, \tag{2.13}$$

there exists an event $\mathcal{A}_*$ with $\mathbb{P}(\mathcal{A}_*) \geq 1 - 2\delta$ such that the iterates generated by Algorithm 1 satisfy the following:

(a) There is a factor $C^{(1)}(d, N, \delta)$ such that

$$\mathbb{E}\left[\tan^2 \angle(\mathbf{u}^{(N)}, \mathbf{u}^*) \, ; \, \mathcal{A}_*\right] \leq C^{(1)}(d, N, \delta) \cdot \frac{\lambda_1}{\lambda_1 - \lambda_2} \sum_{k=2}^{d} \frac{\lambda_k}{\lambda_1 - \lambda_k} \cdot \frac{\log N}{N}. \tag{2.14}$$

The factor $C^{(1)}(d, N, \delta)$ approaches to some absolute positive constant $C_5$ as $d \to \infty, N \to \infty$ when Eqs. (2.12), (2.13) and the following scaling condition hold:

$$d^{-1} \sum_{k=2}^{d} \frac{\lambda_k}{\lambda_1 - \lambda_k} \text{ is bounded away from 0.} \tag{2.15}$$

(b) There is a factor $C^{(2)}(d, N, \delta)$ such that

$$\mathbb{E}\left[\tan^2 \angle(\mathbf{u}^{(N)}, \mathbf{u}^*) ; \mathcal{A}_*\right] \leq C^{(2)}(d, N, \delta) \cdot \frac{\lambda_1 \lambda_2}{(\lambda_1 - \lambda_2)^2} \cdot \frac{(d-1) \log N}{N}. \tag{2.16}$$

The factor $C^{(2)}(d, N, \delta)$ approaches to $C_5$ as $d \to \infty, N \to \infty$ when Eq. (2.12) holds and $\lambda_2/\lambda_1$ is bounded away from 0.

In the above, $b_3, C_5, C_6$ are positive absolute constants.

**Remarks.** Theorem 3 gives an explicit estimate of the convergence rate when a fixed sample size $N$ is known in advance. We note the following:

(i) The conditions (2.12) and (2.13) are satisfied when $d/N^{1-2\varepsilon}$ is sufficiently small. They essentially require that the sample size $N$ be sufficiently large when the dimension $d$ is high.

(ii) The eigenvalue scaling condition Eq. (2.15) requires that the eigenvalues do *not* decay too fast. For example, Eq. (2.15) is satisfied if $\lambda_2/\lambda_1$ is bounded away from 0. For another example, Eq. (2.15) is satisfied if $d^{-1} \sum_{k=1}^{d} \lambda_k/\lambda_1$ is bounded away from 0.

(iii) The choice of stepsize $\beta = \bar{\beta}(N)$ does not involve the dimension $d$. Both constant factors in Eqs. (2.14) and (2.16) satisfy $C^{(i)}(d, N, \delta) \approx C_5, i = 1, 2$ for sufficiently large $d$ and $N$, as long as both the aforementioned scaling conditions hold. From an asymptotic point of view, this implies that as long as $d \to \infty$ and $N \to \infty$ with $d/N^{1-\varepsilon} \to 0$ for some $\varepsilon > 0$, we have $C^{(i)}(d, N, \delta) \to C_5$, and hence $C^{(i)}(d, N, \delta), i = 1, 2$ can be approximately treated as universal constants.

Proofs of Theorems 2 and 3 are deferred to §3.2.

**Matching the Minimax Information Lower Bound.** For an arbitrary constant $\sigma_*^2 > 0$, let $\mathcal{M}(\sigma_*^2, d)$ be the collection of all distributions of $\mathbf{X} = \mathbf{\Sigma}^{1/2} \mathbf{Z}$ satisfying Eq. (2.2) and the following *effective noise variance* condition

$$\frac{\lambda_1 \lambda_2}{(\lambda_1 - \lambda_2)^2} \leq \sigma_*^2; \tag{2.17}$$

see Vu & Lei (2013). For the collection of subgaussian distributions $\mathcal{M}(\sigma_*^2, d)$, Theorem 3.1 of Vu & Lei (2013) established the following *minimax information lower bound* for estimating the corresponding principal component, given by

$$\inf_{\widetilde{\mathbf{u}}^{(N)}} \sup_{\mathbf{X} \in \mathcal{M}(\sigma_*^2, d)} \mathbb{E}\left[\sin^2 \angle(\widetilde{\mathbf{u}}^{(N)}, \mathbf{u}^*)\right] \geq c \cdot \sigma_*^2 \cdot \frac{d-1}{N}, \tag{2.18}$$

for some $c > 0$, where the infimum of $\widetilde{\mathbf{u}}^{(N)}$ is taken over all principal eigenvector estimator using the first $N$ data samples $\boldsymbol{X}^{(1)}, \ldots, \boldsymbol{X}^{(N)}$.

According to Theorem 3(b), we have

$$\sup_{\boldsymbol{X} \in \mathcal{M}(\sigma_*^2, d)} \mathbb{E}\left[\sin^2 \angle(\mathbf{u}^{(N)}, \mathbf{u}^*); \mathcal{A}_*\right] \leq C^{(2)}(d, N, \delta) \cdot \sigma_*^2 \cdot \frac{(d-1)\log N}{N},$$

where $C^{(2)}(d, N, \delta)$ can be approximately treated as a constant. In comparison with Eq. (2.18), Theorem 3 suggests that the online PCA is *nearly optimal*, in the sense that the finite-sample error matches the minimax information lower bound up to a $\log N$ factor with high probability. In addition, online PCA achieves similar accuracy as that of the classical PCA in the *batch* setting (up to $\mathcal{O}(\log N)$). More importantly, the finite-sample bound in Theorem 3 is *nearly global* in the sense that it does *not* require a warm initialization, and hence the near-optimal convergence rate is achieved even using the uniformly random initialization. In summary, our convergence rate result is both *nearly optimal* and *nearly global*.

## 2.3 Comparison with Existing and Concurrent Works

We give a detailed comparison between our finite-sample error bound, i.e., Theorem 3, and related recent results. We emphasize that *all comparable results assume that the samples are uniformly bounded with probability 1*, i.e., $\|\boldsymbol{X}\|^2 \leq B$ for some constant $B > 0$. As a result, their results do not generalize to high-dimensional problems with unbounded distribution. In comparison, we assume that the distribution of $\boldsymbol{X}$ has subgaussian tails, which is much more general than all existing and concurrent works.

A detailed list is given below. For ease of comparison, we rephrase the related results into our setting of uniformly randomized initialization.

(a) The main results of Balsubramani et al. (2013) (Theorem 1.1) can be summarized as follows. Under some technical assumptions, let $\mathbf{u}^{(0)}$ be uniformly sampled from $\mathcal{S}^{d-1}$, let $n_1$ be some starting time, and let the step sizes be $\beta_n = 2(\lambda_1 - \lambda_2)^{-1}(n + n_1)^{-1}$. Then with probability at least 3/4, the iterates satisfy

$$\sin^2 \angle(\mathbf{u}^{(n)}, \mathbf{u}^*) \leq C \cdot \frac{B^2}{(\lambda_1 - \lambda_2)^2} \cdot \frac{1}{n + n_1} + C \cdot d^4 \cdot \left(\frac{n_1}{n + n_1}\right)^2, \qquad (2.19)$$

where $C > 0$ is an absolute constant. This is one of the earliest convergence rate results for the online PCA iteration. Eq. (2.19) does not fully characterize the error's dependence on $\lambda_1, \lambda_2$. It does not match the minimax lower bound (2.18).

(b) The work by De Sa et al. (2015) studies a different but closely related problem on minimizing the spectral error using a stochastic gradient algorithm. The algorithm's *angular part* is equivalent to our online PCA iteration. Their theoretical guarantees are summarized as: Let some technical assumptions hold, let $\mathbf{u}^{(0)}$ be uniformly sampled from $\mathcal{S}^{d-1}$. Given the sample size $N$, by setting $\beta = 16(\lambda_1 - \lambda_2)^{-1} N^{-1}$, the output satisfies with probability at least $3/4$ that

$$\sin^2 \angle(\mathbf{u}^{(N)}, \mathbf{u}^*) \leq C \cdot \frac{B\lambda_1}{(\lambda_1 - \lambda_2)^2} \cdot \frac{d \log N}{N}. \tag{2.20}$$

In addition, De Sa et al. (2015) also proposed an online estimator for the top eigenvalue $\lambda_1$.

(c) Another related work is Shamir (2016a). Shamir (2016a) proves a convergence rate result in terms of the objective values of problem (1.1). We rephrase their main result (Theorem 2) as follows. Let $\mathbf{u}^{(0)}$ be uniformly sampled from $\mathcal{S}^{d-1}$, and let the stepsize be chosen according to the sample size, i.e., $\beta(N) = \lambda_1(\lambda_1 - \lambda_2)^{-1} \cdot N^{-1} \log N$. Then with probability at least $\mathcal{O}(d^{-1})$ the iterates generated by Algorithm 1 satisfy

$$\sin^2 \angle(\mathbf{u}^{(N)}, \mathbf{u}^*) \leq \frac{\mathbf{u}^{*\top} \boldsymbol{\Sigma} \mathbf{u}^* - \mathbf{u}^{(N)\top} \boldsymbol{\Sigma} \mathbf{u}^{(N)}}{\lambda_1 - \lambda_2} \leq C \cdot \frac{B^2}{(\lambda_1 - \lambda_2)^2} \cdot \frac{d \log^2 N}{N}, \tag{2.21}$$

for all $N$ sufficiently large, where $C > 0$ is an absolute constant. The preceding result holds with probability $\mathcal{O}(1/d)$, which approaches to 0 as $d$ grows. In an updated version of Shamir (2016a), the success probability has been improved to $\mathcal{O}(\log^{-1} d)$.

(d) A very recent independent work by Jain et al. (2016) (which is released after the initial submission of the current paper) obtains the following result: Let the $\mathbf{u}^{(0)}$ be uniformly sampled from $\mathcal{S}^{d-1}$, let $n_2$ be some appropriate starting time, and let the stepsize be $\beta_n = (\lambda_1 - \lambda_2)^{-1}(n + n_2)^{-1}$. Then the output $\mathbf{u}^{(n)}$ of Algorithm 1 satisfies with probability at least $3/4$ that

$$\sin^2 \angle(\mathbf{u}^{(n)}, \mathbf{u}^*) \leq C \cdot \frac{B\lambda_1}{(\lambda_1 - \lambda_2)^2} \cdot \frac{1}{n} + C \cdot d \cdot \left(\frac{n_2}{n}\right)^2. \tag{2.22}$$

Here $C > 0$ is an absolute constant.

We summarize all existing rate of convergence results for online PCA in Table 1. In short, the comparable results listed above all require a severely stringent assumption of uniform boundedness. Our results hold under the much more general assumption of subgaussian distributions. We provide the sharpest finite-sample error bound for the online PCA iteration. It is the first result that nearly matches the minimax information lower bound Vu & Lei (2013) up to a polylogarithmic factor of $N$.

## 3 Proofs of Main Results

This section analyzes the convergence of the stochastic iteration generated by Algorithm 1 and provides detailed proofs of the convergence rate results in §2. Arguments for convergence results are separately provided in §3.1 and 3.2, separately for the deterministic initialization case (Theorem 1) and for the uniform randomized initialization case (Theorem 2). §3.2 also provides the argument for the finite-sample result with uniform randomized initialization (Theorem 3). Proofs of all propositions in this section are provided in Section C of Supplementary Material.

| Algorithm | $\sin^2 \angle(\mathbf{u}^{(n)}, \mathbf{u}^*)$ | Optimality |
|---|---|---|
| Matrix Bernstein inequality | $C \cdot \dfrac{B\lambda_1}{(\lambda_1 - \lambda_2)^2} \cdot \dfrac{1}{n}$ | Lower bound |
| Alecton (De Sa et al., 2015) | $C \cdot \dfrac{B\lambda_1 \cdot d}{(\lambda_1 - \lambda_2)^2} \cdot \dfrac{1}{n}$ | No |
| Block power method (Mitliagkas et al., 2013; Hardt & Price, 2014) | $C \cdot \dfrac{B\lambda_1^2}{(\lambda_1 - \lambda_2)^3} \cdot \dfrac{1}{n}$ | No |
| Online PCA, Oja (Balsubramani et al., 2013) | $C \cdot \dfrac{B^2}{(\lambda_1 - \lambda_2)^2} \cdot \dfrac{1}{n}$ | No |
| Online PCA, Oja (Shamir, 2016a) | $C \cdot \dfrac{B^2 \cdot d}{(\lambda_1 - \lambda_2)^2} \cdot \dfrac{1}{n}$ | No |
| Online PCA, Oja (Jain et al., 2016) | $C \cdot \dfrac{B\lambda_1}{(\lambda_1 - \lambda_2)^2} \cdot \dfrac{1}{n}$ | Yes |
| Minimax rate (Vu & Lei, 2013), subgaussian | $C \cdot \dfrac{\lambda_1 \lambda_2 \cdot d}{(\lambda_1 - \lambda_2)^2} \cdot \dfrac{1}{n}$ | Lower bound |
| Online PCA, Oja (this work), subgaussian | $C \cdot \dfrac{\lambda_1}{\lambda_1 - \lambda_2} \sum_{k=2}^{d} \dfrac{\lambda_k}{\lambda_1 - \lambda_k} \cdot \dfrac{1}{n}$ | Yes |

Table 1: Comparable results on the convergence rate of online PCA. Note that our result matches the minimax information lower bound Vu & Lei (2013) in the case where $\lambda_2 = \cdots = \lambda_d$. Our result provides a finer estimate than the minimax lower bound in the more general case where $\lambda_2 \neq \lambda_d$. Note that the constant $C$ hides poly-logarithmic factors of $d$ and $n$.

### 3.1 Proof of Theorem 1

This subsection aims to prove the convergence result for the deterministic initialization case, Theorem 1. To prepare for the proof, first we let the diagonal decomposition of the covariance matrix be

$$\boldsymbol{\Sigma} = \mathbb{E}\left[\boldsymbol{X}\boldsymbol{X}^\top\right] = \mathbf{U}\boldsymbol{\Lambda}\mathbf{U}^\top,$$

where $\boldsymbol{\Lambda} = \text{diag}(\lambda_1, \lambda_2, \ldots, \lambda_d)$ is a diagonal matrix with diagonal entries $\lambda_1, \lambda_2, \ldots, \lambda_d$, and $\mathbf{U}$ is an orthogonal matrix consisting of column eigenvectors of $\boldsymbol{\Sigma}$. Clearly the first column of $\mathbf{U}$ is equal to the principal component $\mathbf{u}^*$. Note that the diagonal decomposition might not be unique, in which case we work with an arbitrary one.

Second, we define the *rescaled samples* as

$$\boldsymbol{Y}^{(n)} = \mathbf{U}^\top \boldsymbol{X}^{(n)}, \qquad \mathbf{v}^{(n)} = \mathbf{U}^\top \mathbf{u}^{(n)}, \qquad \mathbf{v}^* = \mathbf{U}^\top \mathbf{u}^*. \tag{3.1}$$

One can easily verify that

$$\mathbb{E}[\boldsymbol{Y}] = 0, \qquad \mathbb{E}\left[\boldsymbol{Y}\boldsymbol{Y}^\top\right] = \boldsymbol{\Lambda};$$

The principal component of the rescaled random variable $\boldsymbol{Y}$, which we denote by $\mathbf{v}^*$, is equal to $\mathbf{e}_1$, where $\{\mathbf{e}_1, \ldots, \mathbf{e}_d\}$ is the canonical basis of $\mathbb{R}^d$. By applying the linear transformation $\mathbf{U}^\top$ to the stochastic process $\{\mathbf{u}^{(n)}\}$, we obtain an iterative process $\{\mathbf{v}^{(n)} = \mathbf{U}^\top \mathbf{u}^{(n)}\}$ in the rescaled space:

$$\mathbf{v}^{(n)} = \Pi \left\{ \mathbf{v}^{(n-1)} + \beta \boldsymbol{Y}^{(n)} \left(\boldsymbol{Y}^{(n)}\right)^\top \mathbf{v}^{(n-1)} \right\}. \tag{3.2}$$

Moreover, the angle processes associated with $\{\mathbf{u}^{(n)}\}$ and $\{\mathbf{v}^{(n)}\}$ are equivalent, i.e.,

$$\angle(\mathbf{u}^{(n)}, \mathbf{u}^*) = \angle(\mathbf{v}^{(n)}, \mathbf{v}^*).$$

Therefore it would be sufficient to study the rescaled process given by (3.2).

In the proof follows, we use the rescaled samples as in Eq. (3.1) and iterations as in Eq. (3.2). Consider a partition $\mathcal{S}^{d-1} = \mathcal{S}_1 \cup \mathcal{S}_2$ where

$$\mathcal{S}_1 = \left\{\mathbf{v} \in \mathcal{S}^{d-1} : |v_1| < 1/\sqrt{2}\right\}, \qquad \mathcal{S}_2 = \left\{\mathbf{v} \in \mathcal{S}^{d-1} : |v_1| \geq 1/\sqrt{2}\right\}. \tag{3.3}$$

We refer to $\mathcal{S}_1$ and $\mathcal{S}_2$ as the *cold region* and the *warm region*, respectively.

We first focus on Algorithm 1 when the initial estimator lies in the warm region $\mathcal{S}_2$, which we conveniently call *warm start*. Such analysis is crucial in obtaining the correct rate of convergence in the proof of Theorem 1. In terms of the angle $\angle(\mathbf{v}^{(0)}, \mathbf{v}^*)$ this warm start condition is equivalent to

$$\angle(\mathbf{v}^{(0)}, \mathbf{v}^*) \in [0, \pi/4] \cup [3\pi/4, \pi].$$

To avoid uncontrollable variances we need its first coordinate $v_1^{(0)}$ to be bounded away from 0 throughout the algorithm for $N_{\beta,t}$ iterates. We firstly define an auxiliary region

$$\mathcal{S}_3 = \left\{\mathbf{v} \in \mathcal{S}^{d-1} : |v_1| \in [1/3, 1]\right\}, \tag{3.4}$$

and set the stopping time

$$\mathcal{N}_w = \inf\left\{n \geq 0 : \mathbf{v}^{(n)} \in \mathcal{S}_3^c\right\}, \tag{3.5}$$

where $\mathcal{A}^c$ for a generic set $\mathcal{A}$ denotes its complement set. Also, for a positive quantity $M$ to be determined later, let

$$\mathcal{N}_M = \inf\left\{n \geq 1 : \max\left(\max_{1 \leq k \leq d} |Y_k^{(n)}|, |\mathbf{v}^{(n-1)\top}\mathbf{Y}^{(n)}|\right) \geq M^{1/2}\right\}. \tag{3.6}$$

In words, $\mathcal{N}_M$ is the first $n$ such that the maximal absolute coordinate of $\mathbf{Y}^{(n)}$ exceeds $M^{1/2}$, or the inner product of $\mathbf{v}^{(n-1)}$ and $\mathbf{Y}^{(n)}$, in absolute value, exceeds $M^{1/2}$, whichever occurs earlier. It is convenient to define the *rescaled stepsize*

$$\widehat{\beta} = \lambda_1^2(\lambda_1 - \lambda_2)^{-1}\beta. \tag{3.7}$$

We consider the ratio iteration $U_k^{(n)}$ defined as

$$U_k^{(n)} = \frac{v_k^{(n)}}{v_1^{(n)}}. \tag{3.8}$$

Geometrically, we observe that the ratio $U_k^{(n)}$ is the tangent of angle between $\mathbf{v}^{(n)}$ and principal eigenvector $\mathbf{v}^* = \mathbf{e}_1$ after projected onto the two-dimensional subspace spanned by $\mathbf{e}_1$ and $\mathbf{e}_k$ which are the 1st and $k$th canonical unit vectors. Immediately from Eq. (3.8) we have

$$\tan^2 \angle(\mathbf{v}^{(n)}, \mathbf{v}^*) = \sum_{k=2}^{d}\left(U_k^{(n)}\right)^2. \tag{3.9}$$

Our first goal in this subsection is to prove the following

**Proposition 1.** Assume all conditions in Theorem 1 along with the warm start condition $\mathbf{v}^{(0)} \in \mathcal{S}_2$. Then there is a positive constant $b_4$ such that whenever

$$d\widehat{\beta}^{1-2\varepsilon} \leq b_4, \tag{3.10}$$

for any fixed $t > 0$ and for $k = 2, \ldots, d$ we have with probability $\geq 1 - 2N_{\beta,t}^* \exp\left(-C_{1,P}\widehat{\beta}^{-2\varepsilon}\right)$, either $(\mathcal{N}_M \leq N_{\beta,t}^*)$ occurs, or

$$\sup_{n \leq N_{\beta,t}^* \wedge \mathcal{N}_w} \left| U_k^{(n)} - U_k^{(0)} \left(1 - \beta(\lambda_1 - \lambda_k)\right)^n \right| \leq C'_{1,P}\widehat{\beta}^{0.5-\varepsilon}.$$

Here $C_{1,P}$ and $C'_{1,P}$ are positive constants.

Proof of Proposition 1 is deferred to §C.1. Proposition 1 shows that from warm start, for each coordinate $k = 2, \ldots, d$, the $U_k^{(n)}$ approximately decays geometrically at rate $1 - \beta(\lambda_1 - \lambda_k)$, or *bad event* $(\mathcal{N}_M \leq N_{\beta,t}^*)$ occurs.

The next lemma controls the occurrence of the bad event by pinning down the choice of $M$.

**Lemma 1.** Assume all conditions in Theorem 1. Suppose we choose $M$ as

$$M \equiv \lambda_1 \widehat{\beta}^{-2\varepsilon}. \tag{3.11}$$

Then for each fixed $N \geq 1$

$$\mathbb{P}\left(\mathcal{N}_M \leq N\right) \leq 2(d+1)N \exp\left(-\widehat{\beta}^{-2\varepsilon}\right). \tag{3.12}$$

The proof of Lemma 1 is provided in §D.1. Using Proposition 1 and Lemma 1, one can obtain the rate of convergence under more careful second moment estimates, as follows.

**Proposition 2.** Assume all conditions in Theorem 1 along with the warm start condition $\mathbf{v}^{(0)} \in \mathcal{S}_2$. Then there is a positive constant $b_4$ such that whenever Eq. (3.10) is satisfied, there exists a high-probability event $\mathcal{H}_0$ satisfying

$$\mathbb{P}(\mathcal{H}_0) \geq 1 - 4dN_{\beta,t}^* \exp\left(-C_{2,P}\widehat{\beta}^{-2\varepsilon}\right), \tag{3.13}$$

such that for fixed $k = 2, \ldots, d$ and $n \in [N_{\beta,1}^*, N_{\beta,t}^*]$,

$$\mathbb{E}\left[\left(U_k^{(n)}\right)^2 ; \mathcal{H}_0\right] \leq (1 - \beta(\lambda_1 - \lambda_k))^{2n}\left(U_k^{(0)}\right)^2 + \frac{C'_{2,P}(\lambda_1\lambda_k + \lambda_1^2 \cdot \widehat{\beta}^{0.5-4\varepsilon})}{\lambda_1 - \lambda_k} \cdot \beta \tag{3.14}$$

Here $C_{2,P}$ and $C'_{2,P}$ are positive constants.

Proof of Proposition 2 is deferred to §C.2. Note from Eq. (3.13) $\mathcal{H}_0$ occurs with high probability as $\beta \to 0$. With $\widehat{\beta}$ defined in Eq. (3.7) the second term on the right hand of Eq. (3.14) matches the second moment estimates as in Eq. (2.7).

Second, we try to relax the warm start condition. Define the stopping time

$$\mathcal{N}_c = \inf\left\{n \geq 0 : \mathbf{v}^{(n)} \in \mathcal{S}_2\right\}. \tag{3.15}$$

In words, $\mathcal{N}_c$ is the first $n$ such that the iterate enters the warm region. In the case of warm start, $\mathcal{N}_c = 0$. We have the following Proposition 3 which upper bound the stopping time $\mathcal{N}_c$, whose proof is deferred to §C.3.

**Proposition 3.** Assume all conditions in Theorem 1. There exist a positive constant $b_5 < \ln^2 2/16$ such that if

$$d\widehat{\beta}^{1-2\varepsilon} \leq b_5/c^*, \tag{3.16}$$

then we have

$$\mathbb{P}\left(\mathcal{N}_c \leq N_\beta^o(c^*)\right) \geq 1 - 2(d+2)N_\beta^o(c^*)\exp\left(-C_{3,P}\widehat{\beta}^{-2\varepsilon}\right), \tag{3.17}$$

where $N_\beta^o(c^*)$ is defined as in Eq. (2.4).

We are now ready prove Theorem 1. Note that under the setting of Theorem 1, Algorithm 1 starts from $\mathbf{v}^{(0)} \in \mathcal{S}^{d-1}$ where $\tan^2\angle(\mathbf{v}^{(0)}, \mathbf{v}^*) \leq c^*d$. Running the algorithm for $N_\beta^o(c^*) \wedge \mathcal{N}_c$ steps we know from Proposition 3 that the iterate $\mathbf{v}^{(N_\beta^o(c^*)\wedge\mathcal{N}_c)}$ lies in $\mathcal{S}_2$ with high probability. By strong Markov property, the iterates generated by Algorithm 1 have the same law as the one starting from $\mathbf{v}^{(N_\beta^o(c^*)\wedge\mathcal{N}_c)}$, and hence Proposition 2 can be applied. We follow this reasoning and detail the proof as in below.

*Proof of Theorem 1.* Let the events

$$\begin{aligned}
\mathcal{H}_{*,1} &= \left(\mathcal{N}_c \leq N_\beta^o(c^*)\right), \\
\mathcal{H}_{*,2} &= \left(\sup_{n \in [N_{\beta,1}^*, N_{\beta,t}^*]} \left|U_k^{(n+\mathcal{N}_c)}\right| \leq 2\widehat{\beta}^{0.5-\varepsilon} \text{ for all } k = 2, \ldots, d\right), \\
\mathcal{H}_* &= \mathcal{H}_{*,1} \cap \mathcal{H}_{*,2}.
\end{aligned} \tag{3.18}$$

On $\mathcal{H}_* \cap (\mathcal{N}_c = n_o)$ where $n_o \leq N_\beta^o(c^*)$ is a fixed time, from the definition in Eq. (3.15) we have $\mathbf{v}^{(\mathcal{N}_c \wedge N_\beta^o(c^*))} = \mathbf{v}^{(n_o)} \in \mathcal{S}_2$ and hence from Eq. (3.9)

$$\tan^2\angle(\mathbf{v}^{(n_o)}, \mathbf{v}^*) = \sum_{k=2}^d \left(U_k^{(n_o)}\right)^2 \leq 1. \tag{3.19}$$

Therefore for

$$n \in [N_{\beta,1}^* + N_\beta^o(c^*), N_{\beta,t}^*] \subseteq [N_{\beta,1}^* + n_o, N_{\beta,t}^* + n_o],$$

and $\beta$ satisfying both Eqs. (3.10) and (3.16), we utilize the Markov property and Proposition 2 and conclude that for $\beta$ satisfying Eq. (2.5) where $b_1 = b_4 \wedge b_5$

$$\mathbb{E}\left[\tan^2 \angle(\mathbf{v}^{(n)}, \mathbf{v}^*)\mathbf{1}_{\mathcal{H}_*} \,|\, \mathcal{N}_c = n_o\right] = \sum_{k=2}^{d} \mathbb{E}\left[\left(U_k^{(n)}\right)^2 \mathbf{1}_{\mathcal{H}_*} \,\bigg|\, \mathcal{N}_c = n_o\right]$$

$$\leq \sum_{k=2}^{d} (1-\beta(\lambda_1 - \lambda_k))^{2(n-n_o)} \left(U_k^{(n_o)}\right)^2 + \sum_{k=2}^{d} \frac{C'_{2,P}(\lambda_1\lambda_k + \lambda_1^2 \cdot \widehat{\beta}^{0.5-4\varepsilon})}{2(\lambda_1 - \lambda_k)}\beta$$

$$\leq (1-\beta(\lambda_1 - \lambda_2))^{2(n-N_\beta^o(c^*))} + \sum_{k=2}^{d} \frac{C'_{2,P}(\lambda_1\lambda_k + \lambda_1^2 \cdot \widehat{\beta}^{0.5-4\varepsilon})}{2(\lambda_1 - \lambda_k)}\beta.$$

Eqs. (3.9), (3.19) have been applied in the above estimate. The bound for

$$\mathbb{E}\left[\tan^2 \angle(\mathbf{v}^{(n)}, \mathbf{v}^*)\mathbf{1}_{\mathcal{H}_*} \,|\, \mathcal{N}_c\right]$$

was argued for $\mathcal{N}_c = n_o \leq N_\beta^o(c^*)$ but is trivially valid on $(\mathcal{N}_c > N_\beta^o(c^*))$, since $\mathbf{1}_{H_*} = 0$ on the latter event. Furthermore, the bound on the right hand is deterministic and hence independent of $\mathcal{N}_c$. Taking expectation again gives

$$\mathbb{E}\left[\tan^2 \angle(\mathbf{v}^{(n)}, \mathbf{v}^*) \,;\, \mathcal{H}_*\right] \leq (1-\beta(\lambda_1 - \lambda_2))^{2(n-N_\beta^o(c^*))}$$
$$+ \sum_{k=2}^{d} \frac{C'_{2,P}(\lambda_1\lambda_k + \lambda_1^2 \cdot \widehat{\beta}^{0.5-4\varepsilon})}{2(\lambda_1 - \lambda_k)}\beta.$$

Therefore Eq. (2.7) establishes with $C_2 = C'_{2,P}$.

To ensure that the event $\mathcal{H}_*$ defined in Eq. (3.18) satisfies the probability estimate in Eq. (2.6), note Eqs. (3.13) and (3.17) together imply

$$\mathbb{P}(\mathcal{H}_*) \geq \left[1 - 2(d+2)N_\beta^o(c^*)\exp\left(-C_{3,P}\widehat{\beta}^{-2\varepsilon}\right)\right]\left[1 - 4dN^*_{\beta,t}\exp\left(-C_{2,P}\widehat{\beta}^{-2\varepsilon}\right)\right]$$
$$\geq 1 - 2(d+2)N_\beta^o(c^*)\exp\left(-C_{3,P}\widehat{\beta}^{-2\varepsilon}\right) - 4dN^*_{\beta,t}\exp\left(-C_{2,P}\widehat{\beta}^{-2\varepsilon}\right),$$

Letting $C_0 = C_{3,P}$ and $C_1 = C_{2,P}$ this implies that Eq. (2.6) holds. Proof of Theorem 1 is accomplished.

$\square$

### 3.2 Proofs of Theorem 2 and 3

Next we prove the main results of the uniform randomized initialization case, Theorems 2 and 3. The following lemma basically says that the initial estimator condition in Theorem 1 is satisfied with high probability, where the coefficient $c^* = C^*\delta^{-2}$ is quadratically inverse proportional to the error probability $\delta$.

**Lemma 2.** Given any $\delta > 0$, if $\mathbf{u}^{(0)}$ is sampled uniformly at random from $\mathcal{S}^{d-1}$ in $\mathbb{R}^d$ then there exists a constant $C^* > 1$ independent of $\delta$ and $d$ such that

$$\mathbb{P}\left(\tan^2 \angle(\mathbf{u}^{(0)}, \mathbf{u}^*) \leq C^*\delta^{-2}d\right) \geq 1 - \delta.$$

Proof of Lemma 2 is provided in §D.2. We are ready to prove Theorem 2.

*Proof of Theorem 2.* Let the event

$$\mathcal{A}'_* = \left(\tan^2 \angle(\mathbf{u}^{(0)}, \mathbf{u}^*) \leq C^* \delta^{-2} d\right),$$

and recall $\widehat{\beta}$ is as defined in Eq. (3.7). Since the initial estimator $\mathbf{u}^{(0)}$ is sampled uniformly at random from $\mathcal{S}^{d-1}$, Lemma 2 indicates $\mathbb{P}(\mathcal{A}'_*) \geq 1 - \delta$, and our approach is to apply Theorem 1 with $c^* = C^* \delta^{-2}$.

Firstly, Eq. (3.16) in Proposition 3, one has $d\widehat{\beta}^{1-2\varepsilon} \leq (4C^* \delta^{-2})^{-1}$ i.e.

$$\log\left(4C^* \delta^{-2} d\right) \leq (1 - 2\varepsilon) \log(\widehat{\beta}^{-1}).$$

Combined with the definitions in Eqs. (2.3) and (2.4), this implies

$$N_\beta^o(C^* \delta^{-2}) \leq N_{\beta, 1-2\varepsilon}^* \leq N_{\beta, 1}^* - 1. \tag{3.20}$$

Conditioning on $\mathcal{A}'_*$ which satisfies the initial condition in Theorem 1, Eq. (3.20) guarantees the existence of the event $\mathcal{H}_*$ such that

$$\begin{aligned}
\mathbb{P}(\mathcal{H}_* \mid \mathcal{A}'_*) &\geq 1 - 2(d+2)N_\beta^o(c^*) \exp\left(-C_0 \widehat{\beta}^{-2\varepsilon}\right) - 4dN_{\beta,t}^* \exp\left(-C_1 \widehat{\beta}^{-2\varepsilon}\right) \\
&\geq 1 - 4dN_{\beta,t+1}^* \exp\left(-C_3 \widehat{\beta}^{-2\varepsilon}\right) \geq 1 - \delta,
\end{aligned} \tag{3.21}$$

where $C_3 \equiv C_0 \wedge C_1$. Thus immediately

$$\mathbb{P}(\mathcal{A}'_* \cap \mathcal{H}_*) = \mathbb{P}(\mathcal{A}'_*)\mathbb{P}(\mathcal{H}_* \mid \mathcal{A}'_*) \geq (1 - \delta)^2 \geq 1 - 2\delta.$$

Finally, for all $\beta$ satisfying $d\widehat{\beta}^{1-2\varepsilon} \leq (b_1/C^*)\delta^2$,

$$\begin{aligned}
\mathbb{E}\left[\tan^2 \angle(\mathbf{u}^{(n)}, \mathbf{u}^*); \mathcal{A}'_* \cap \mathcal{H}_*\right] &\leq \mathbb{E}\left[\tan^2 \angle(\mathbf{u}^{(n)}, \mathbf{u}^*)\mathbf{1}_{\mathcal{H}_*} \mid \mathcal{A}'_*\right] \\
&\leq 4C^* \delta^{-2} d \left(1 - \beta(\lambda_1 - \lambda_2)\right)^{2n} + C_2 \sum_{k=2}^{d} \frac{\lambda_1 \lambda_k + \lambda_1^2 \widehat{\beta}^{0.5-4\varepsilon}}{\lambda_1 - \lambda_k} \cdot \beta.
\end{aligned}$$

where we applied from Eq. (2.4) that $(1 - \beta(\lambda_1 - \lambda_2))^{-2N_\beta^o(C^* \delta^{-2})} \leq \left(4C^* \delta^{-2} d\right)^2$.

To summarize all above, setting $C_4 = (4C^*)^2 \vee C_2$, $b_2 = [b_1 \wedge (1/4)]/C^*$ and $\mathcal{A}_* = \mathcal{A}'_* \cap \mathcal{H}_*$ concludes Theorem 2.

$\square$

We finalize this paper by proving the finite-sample result, Theorem 3.

*Proof of Theorem 3.* As Eq. (2.12) holds, choosing $\beta$ as in Eq. (2.11), then if for some positive constants $b_2 > 0$ and $\varepsilon = \varepsilon_0 \in (0, 1/8)$, $\beta$ and $\delta > 0$ satisfies Eq. (2.8) which is translated to

$$d \leq b_2 \delta^2 \left[\frac{\lambda_1^2}{(\lambda_1 - \lambda_2)^2} \cdot \frac{\log N}{N}\right]^{2\varepsilon - 1}. \tag{3.22}$$

Hence Theorem 2 applies, and it is not hard to verify when $\beta = \bar{\beta}(N)$ defined in Eq. (2.11) one has immediately $N \in [N^*_{\beta,1}, N^*_{\beta,2}]$, and hence there exists an event $\mathcal{A}_*$ with $\mathbb{P}(\mathcal{A}_*) \geq 1 - 2\delta$, such that the iterates generated by Algorithm 1 satisfy Eq. (2.10) for $n \in [N^*_{\beta,2-\varepsilon}, N^*_{\beta,t}]$. We consider the two regimes separately.

(a) Plugging in $\beta = \bar{\beta}(N)$, and hence to obtain Eq. (2.14)

$$\mathbb{E}\left[\tan^2 \angle(\mathbf{u}^{(N)}, \mathbf{u}^*) ; \mathcal{A}_*\right] \leq C^o(d, N, \delta) \cdot \frac{\lambda_1}{\lambda_1 - \lambda_2} \sum_{k=2}^{d} \frac{\lambda_k}{\lambda_1 - \lambda_k} \cdot \frac{\log N}{N},$$

where the factor $C^o(d, N, \delta)$ takes the form

$$C^o(d, N, \delta) = \left(C_4 \delta^{-2} d \left(1 - \bar{\beta}(N)(\lambda_1 - \lambda_2)\right)^{2N}\right.$$
$$\left. + C_4 \sum_{k=2}^{d} \frac{\lambda_1 \lambda_k + \lambda_1^2 \left[\lambda_1^2 (\lambda_1 - \lambda_2)^{-1} \bar{\beta}(N)\right]^{0.5 - 4\varepsilon}}{(\lambda_1 - \lambda_k)(\lambda_1 - \lambda_2)} \cdot \frac{\log N}{N}\right) \quad (3.23)$$
$$\cdot \left(\frac{\lambda_1}{\lambda_1 - \lambda_2} \sum_{k=2}^{d} \frac{\lambda_k}{\lambda_1 - \lambda_k} \cdot \frac{\log N}{N}\right)^{-1} = \mathrm{I} + \mathrm{II}.$$

From Eqs. (3.22) and (2.11)

$$\mathrm{I} \equiv C_4 \delta^{-2} d \left(1 - \bar{\beta}(N)(\lambda_1 - \lambda_2)\right)^{2N} \cdot \left(\frac{\lambda_1}{\lambda_1 - \lambda_2} \sum_{k=2}^{d} \frac{\lambda_k}{\lambda_1 - \lambda_k} \cdot \frac{\log N}{N}\right)^{-1}$$
$$\leq C_4 \delta^{-2} d \exp\left(-2N \cdot \frac{\log N}{N}\right) \cdot \left(\frac{\lambda_1 \lambda_2}{(\lambda_1 - \lambda_2)^2} \cdot \frac{\log N}{N}\right)^{-1}$$
$$\leq C_4 b_2 \left[\frac{\lambda_1^2}{(\lambda_1 - \lambda_2)^2} \cdot \frac{\log N}{N}\right]^{2\varepsilon - 1} N^{-2} \cdot \left(\frac{\log N}{N}\right)^{-1} \cdot \left(\frac{\lambda_1 \lambda_2}{(\lambda_1 - \lambda_2)^2}\right)^{-1}$$
$$\leq C_4 b_2 \left(\frac{\lambda_1}{\lambda_1 - \lambda_2}\right)^{4\varepsilon - 2} \left(\frac{\lambda_1 \lambda_2}{(\lambda_1 - \lambda_2)^2}\right)^{-1} \frac{1}{N^{2\varepsilon} \log^{2 - 2\varepsilon} N},$$

and also

$$\mathrm{II} \equiv C_4 \sum_{k=2}^{d} \frac{\lambda_1 \lambda_k + \lambda_1^2 \left[\lambda_1^2 (\lambda_1 - \lambda_2)^{-1} \bar{\beta}(N)\right]^{0.5 - 4\varepsilon}}{(\lambda_1 - \lambda_k)(\lambda_1 - \lambda_2)} \cdot \left(\frac{\lambda_1}{\lambda_1 - \lambda_2} \sum_{k=2}^{d} \frac{\lambda_k}{\lambda_1 - \lambda_k}\right)^{-1}$$
$$= C_4 \left(1 + \left[\lambda_1^2 (\lambda_1 - \lambda_2)^{-1} \bar{\beta}(N)\right]^{0.5 - 4\varepsilon} \left[d \left(\sum_{k=2}^{d} \frac{\lambda_k}{\lambda_1 - \lambda_k}\right)^{-1} + 1\right]\right).$$

Both terms on the RHS of Eq. (3.24) have been estimated, and we let

$$C^{(1)}(d, N, \delta) = C_4 b_2 \left(\frac{\lambda_1}{\lambda_1 - \lambda_2}\right)^{4\varepsilon - 2} \left(\frac{\lambda_1 \lambda_2}{(\lambda_1 - \lambda_2)^2}\right)^{-1} \frac{1}{N^{2\varepsilon} \log^{2 - 2\varepsilon} N}$$
$$+ C_4 + C_4 \left[\lambda_1^2 (\lambda_1 - \lambda_2)^{-1} \bar{\beta}(N)\right]^{0.5 - 4\varepsilon} \left[d \left(\sum_{k=2}^{d} \frac{\lambda_k}{\lambda_1 - \lambda_k}\right)^{-1} + 1\right],$$

then we conclude immediately that as $d, N \to \infty$ when both Eqs. (2.12) and (2.15) hold, the term goes to $C_4$. Letting $C_5 \equiv C_4$ completes the proof.

(b) For Eq. (2.16) we note $\lambda_k \leq \lambda_2$ for $k \geq 2$ and

$$\mathbb{E}\left[\tan^2 \angle(\mathbf{u}^{(N)}, \mathbf{u}^*); \mathcal{A}_*\right] \leq C^\omega(d, N, \delta) \cdot \frac{\lambda_1 \lambda_2}{(\lambda_1 - \lambda_2)^2} \cdot \frac{(d-1)\log N}{N},$$

where the factor $C^\omega(d, N, \delta)$ has, via Eqs. (3.22) and (2.11),

$$\begin{aligned}
C^\omega(d, N, \delta) &= \left(C_4 \delta^{-2} d \left(1 - \bar{\beta}(N)(\lambda_1 - \lambda_2)\right)^{2N}\right. \\
&\quad \left. + C_4 \frac{\lambda_1 \lambda_2 + \lambda_1^2 \left[\lambda_1^2(\lambda_1 - \lambda_2)^{-1} \bar{\beta}(N)\right]^{0.5 - 4\varepsilon}}{(\lambda_1 - \lambda_2)^2} \cdot \frac{(d-1)\log N}{N}\right) \\
&\quad \cdot \left(\frac{\lambda_1 \lambda_2}{(\lambda_1 - \lambda_2)^2} \cdot \frac{(d-1)\log N}{N}\right)^{-1} \\
&\leq 2 C_4 \delta^{-2} \exp\left(-2N \cdot \frac{\log N}{N}\right) \cdot \left(\frac{\lambda_1 \lambda_2}{(\lambda_1 - \lambda_2)^2} \cdot \frac{\log N}{N}\right)^{-1} \\
&\quad + C_4 \left(1 + \left[\lambda_1^2(\lambda_1 - \lambda_2)^{-1} \bar{\beta}(N)\right]^{0.5 - 4\varepsilon} \cdot \frac{\lambda_1}{\lambda_2}\right) \\
&\leq 2 C_4 \delta^{-2} \left(\frac{\lambda_1 \lambda_2}{(\lambda_1 - \lambda_2)^2}\right)^{-1} \frac{1}{N \log N} + C_4 + C_4 \left[\lambda_1^2(\lambda_1 - \lambda_2)^{-1} \bar{\beta}(N)\right]^{0.5 - 4\varepsilon} \cdot \frac{\lambda_1}{\lambda_2}.
\end{aligned}$$
(3.24)

Letting $C^{(2)}(d, N, \delta)$ be the last line above concludes that as $d, N \to \infty$ when Eq. (2.12) holds and $\lambda_2/\lambda_1$ being bounded away from 0, the term goes to $C_4 = C_5$. This completes our proof.

$\square$

## 4 Summary

In this paper, we provide an explicit convergence rate analysis for the online PCA iteration with subgaussian samples. Our convergence rate result is *nearly optimal* in the sense that the output can be viewed as a running estimator that nearly attains the minimax information lower bound for subgaussian PCA. Furthermore, our convergence rate result is *nearly global*, in the sense that finite-sample error bound holds with high probability even if the initial solution is uniformly sampled from the unit sphere. One direction for future research is to develop and analyze online PCA method to estimate the top $k$ eigenvectors. Another direction for research is to extend the analysis beyond PCA to a broader class of low-rank estimation problems.

**Acknowledgements** The authors would like to thank Ruoyu Sun and Zhaoran Wang for inspiring discussions on this work. The authors would also like to thank the two anonymous referees for their insightful comments.

# A  Analysis of Algorithm Increments

Throughout the Supplementary Material we use the following notations: (i) The $C$'s with subscripts denotes some positive constants; (ii) The $C$, $C'$, $C''$'s (*without* subscripts) are positive constants whose values may change between lines; (iii) The $\mathbf{v} \equiv \mathbf{v}^{(n-1)}$, $\widehat{\mathbf{v}} \equiv \mathbf{v}^{(n)}$ and $\boldsymbol{Y} \equiv \boldsymbol{Y}^{(n)}$, and the corresponding coordinates are in the same fashion.

To analyze the algorithm we need understand the increments on each coordinate at each step.

**Proposition 4.** Under Assumption 1, for each $k = 1, 2, \ldots, d$ and $n \geq 0$ we have for all $dM\beta \leq 1/3$ the following:

(i) There exists a random variable $Q_{k,n}$ such that on the event $(\mathcal{N}_M > n)$, $|Q_{k,n}| \leq C_{4,1}M^2\beta^2$ almost surely, and the increment at time $n$ on coordinate $k$ that $v_k^{(n)} - v_k^{(n-1)}$ can be expressed as

$$v_k^{(n)} - v_k^{(n-1)} = \beta\left[(\mathbf{v}^\top \boldsymbol{Y})Y_k - v_k(\mathbf{v}^\top \boldsymbol{Y})^2\left(1 + \frac{\beta}{2}\|\boldsymbol{Y}\|^2\right)\right] + Q_{k,n}, \tag{A.1}$$

where we denote $\mathbf{v} = \mathbf{v}^{(n-1)}$ and $\boldsymbol{Y} = \boldsymbol{Y}^{(n)}$.

(ii) On $(\mathcal{N}_M > n)$ the increment has the following bound

$$\left|v_k^{(n)} - v_k^{(n-1)}\right| \leq C_{4,2}M\beta. \tag{A.2}$$

## A.1  Proof of Proposition 4

This subsection is devoted to the proof of Proposition 4. We first come to show

**Lemma 3.** For each $n \geq 0$, on the event $(\mathcal{N}_M > n)$

$$\left|\|\mathbf{v} + \beta(\mathbf{v}^\top \boldsymbol{Y})\boldsymbol{Y}\|^{-1} - 1 + \beta(\mathbf{v}^\top \boldsymbol{Y})^2 + \frac{1}{2}\beta^2(\mathbf{v}^\top \boldsymbol{Y})^2\|\boldsymbol{Y}\|^2\right| \leq C_3\beta^2(\mathbf{v}^\top \boldsymbol{Y})^4.$$

*Proof.* Since

$$\|\mathbf{v} + \beta(\mathbf{v}^\top \boldsymbol{Y})\boldsymbol{Y}\|^{-1} = \left(1 + 2\beta(\mathbf{v}^\top \boldsymbol{Y})^2 + \beta^2(\mathbf{v}^\top \boldsymbol{Y})^2\|\boldsymbol{Y}\|^2\right)^{-1/2}, \tag{A.3}$$

Taylor expansion implies for $|x| < 1$

$$(1+x)^{-1/2} = \sum_{n=0}^{\infty} \binom{-1/2}{n} x^n = 1 - \frac{1}{2}x + \frac{3}{8}x^2 - \frac{5}{16}x^3 + \cdots,$$

which is an alternating series for $x \in [0, 1)$, whereas the absolute terms approach to 0 monotonically

$$\left|\binom{-1/2}{n+1} x^{n+1}\right| \leq \left|\binom{-1/2}{n} x^n\right|.$$

Hence the error bound gives for all $x \in [0, 1)$

$$\left|(1+x)^{-1/2} - 1 + \frac{1}{2}x\right| \leq \frac{3}{8}x^2. \tag{A.4}$$

Noting $|\mathbf{v}^\top \boldsymbol{Y}| \leq \|\boldsymbol{Y}\|$ we have for all $\beta$

$$2\beta(\mathbf{v}^\top \boldsymbol{Y})^2 + \beta^2(\mathbf{v}^\top \boldsymbol{Y})^2 \|\boldsymbol{Y}\|^2 \leq 2M\beta + M^2\beta^2.$$

The above display is strictly less than 1 when $\beta \leq (3M)^{-1}$, and hence (A.4) applies. Combined with (A.3) we have

$$\left| \|\mathbf{v} + \beta(\mathbf{v}^\top \boldsymbol{Y})\boldsymbol{Y}\|^{-1} - 1 + \frac{1}{2}\left(2\beta(\mathbf{v}^\top \boldsymbol{Y})^2 + \beta^2(\mathbf{v}^\top \boldsymbol{Y})^2\|\boldsymbol{Y}\|^2\right)\right| \leq \frac{3}{8}\left(3\beta(\mathbf{v}^\top \boldsymbol{Y})^2\right)^2.$$

Noticing $|v^\top \boldsymbol{Y}| \leq \|\boldsymbol{Y}\|$, triangle inequality suggests

$$\left| \|\mathbf{v} + \beta(\mathbf{v}^\top \boldsymbol{Y})\boldsymbol{Y}\|^{-1} - 1 + \beta(\mathbf{v}^\top \boldsymbol{Y})^2 \right| \leq C\beta^2 \|\boldsymbol{Y}\|^4 \leq CM^2\beta^2,$$

completing the proof.

$\square$

*Proof of Proposition 4.* Setting $Q = \|\mathbf{v} + \beta(\mathbf{v}^\top \boldsymbol{Y})\boldsymbol{Y}\|^{-1} - 1 + \beta(\mathbf{v}^\top \boldsymbol{Y})^2$. Then

$$\begin{aligned}
\widehat{v}_k - v_k &= \|\mathbf{v} + \beta(\mathbf{v}^\top \boldsymbol{Y})\boldsymbol{Y}\|^{-1}\left(v_k + \beta \mathbf{v}^\top \boldsymbol{Y} Y_k\right) - v_k \\
&= \left(1 - \beta(\mathbf{v}^\top \boldsymbol{Y})^2 + Q\right)\left(v_k + \beta \mathbf{v}^\top \boldsymbol{Y} Y_k\right) - v_k \\
&= \beta\left((\mathbf{v}^\top \boldsymbol{Y}) Y_k - v_k(\mathbf{v}^\top \boldsymbol{Y})^2\right) + Q_{k,n},
\end{aligned}$$

where

$$Q_{k,n} = \left(v_k + \beta \mathbf{v}^\top \boldsymbol{Y} Y_k\right)Q - \beta^2(\mathbf{v}^\top \boldsymbol{Y})^3 Y_k. \tag{A.5}$$

Note the term

$$\beta\left[(\mathbf{v}^\top \boldsymbol{Y})Y_k - v_k(\mathbf{v}^\top \boldsymbol{Y})^2\right]$$

is absolutely bounded by $2M\beta$, and taking expectation gives

$$\begin{aligned}
\mathbb{E}\left[(\mathbf{v}^\top \boldsymbol{Y})Y_k - v_k(\mathbf{v}^\top \boldsymbol{Y})^2\right] &= v_k \lambda_k - v_k \mathbb{E}(\mathbf{v}^\top \boldsymbol{Y})^2 \\
&= v_k \lambda_k - v_k \mathbf{v}^\top \mathbb{E}(\boldsymbol{Y}\boldsymbol{Y}^\top)\mathbf{v}^\top = v_k \left(\lambda_k - \mathbf{v}^\top \boldsymbol{\Lambda} \mathbf{v}\right).
\end{aligned}$$

To this stage, we have verified

$$\widehat{v}_k - v_k = \beta\left((\mathbf{v}^\top \boldsymbol{Y})Y_k - v_k(\mathbf{v}^\top \boldsymbol{Y})^2\right) + Q_{k,n}. \tag{A.6}$$

(A.1) along with (A.2) in Proposition 4 can be concluded if

$$|Q_{k,n}| \leq CM^2\beta^2. \tag{A.7}$$

To conclude (A.7), note that $\beta \leq (3M)^{-1}$ and hence

$$\left|v_k + \beta \mathbf{v}^\top \boldsymbol{Y} Y_k\right| \leq 1 + \beta M \leq \frac{4}{3}.$$

Lemma 3 implies
$$|Q| \leq C_3\beta^2(\mathbf{v}^\top \mathbf{Y})^4 \leq C_3 M^2 \beta^2.$$

Therefore the first term on the right hand of (A.5) is absolutely bounded by $2C_3 M^2 \beta^2$. For the second term in (A.5) we have
$$|\beta^2(\mathbf{v}^\top \mathbf{Y})^3 Y_k| \leq M^2 \beta^2.$$

We have verified (A.7) by taking $C = 2C_3 + 1$, which completes all the proof of Proposition 4.

□

## B  Subgaussian Bound

For further analysis, we provide in this section the concepts and basic properties of subgaussian and subexponential random variables in §B.2, and finally prove the concentration inequality Lemma 8 in §B.3.

### B.1  Subgaussian random variables

Recall that in Definition 1 we have for subgaussian random variable $Y$ that for each $p \geq 1$
$$(\mathbb{E}|Y|^p)^{1/p} \leq p^{1/2}\|Y\|_{\psi_2}. \tag{B.1}$$

We conclude the following subgaussian-tail property.

**Lemma 4** (Property of subgaussian distributions)**.** Let $Y$ be a random variable satisfying (B.1). We have for all $t \geq 0$
$$\mathbb{P}(|Y| \geq t) \leq 2\exp\left(-t^2/\|Y\|_{\psi_2}^2\right). \tag{B.2}$$

*Proof.* Exponentiating and using Markov's inequality as well as (B.1), we have
$$\mathbb{P}(|Y| \geq t) = \mathbb{P}\left(\exp\left(Y^2/\|Y\|_{\psi_2}^2\right) \geq \exp\left(t^2/\|Y\|_{\psi_2}^2\right)\right)$$
$$\leq \exp\left(-t^2/\|Y\|_{\psi_2}^2\right) \mathbb{E}\exp\left(Y^2/\|Y\|_{\psi_2}^2\right) \leq 2\exp\left(-t^2/\|Y\|_{\psi_2}^2\right).$$

This proves the subgaussian tail property (B.2).

□

Based on the above definition we first characterize the fourth moments of subgaussian variables using the following lemma.

**Lemma 5.** Let $\mathbf{Z}$ be vector-subgaussian with $\sup_{\mathbf{v}\in\mathcal{S}^{d-1}} \|\mathbf{v}^\top \mathbf{Z}\|_{\psi_2} \leq 1$. We conclude for any vector $\mathbf{w}$ that
$$\mathbb{E}\left(\mathbf{w}^\top \mathbf{Z}\right)^4 \leq 16\|\mathbf{w}\|^4. \tag{B.3}$$

*Proof.* Taking $p = 4$ in (B.1), we have from $\sup_{\mathbf{v}\in\mathcal{S}^{d-1}} \|\mathbf{v}^\top \mathbf{Z}\|_{\psi_2} \leq 1$ that
$$\|\mathbf{w}\|^{-4} \cdot \mathbb{E}\left(\mathbf{w}^\top \mathbf{Z}\right)^4 \leq 4^2 \left\|\frac{\mathbf{w}^\top}{\|\mathbf{w}\|} \cdot \mathbf{Z}\right\|_{\psi_2}^4 \leq 16.$$

This concludes the proof.

□

## B.2 Subexponential random variables

In this subsection we introduce the concept of subexponential norm as well as random variables that have subexponential tails. We follow Vershynin (2010) and define the $\psi_1$-norm as follows.

**Definition 3** (Subexponential random variables). The subexponential norm of $X$ is defined as

$$\|X\|_{\psi_1} \equiv \sup_{1 \leq p < \infty} p^{-1} \left(\mathbb{E}|X|^p\right)^{1/p}.$$

For $X$ that satisfies $\|X\|_{\psi_1} < \infty$ is called a *subexponential* random variable.

**Remark.** It is straightforward to verify via Minkowski's inequality the following triangle inequality

$$\|X + Y\|_{\psi_1} \leq \|X\|_{\psi_1} + \|Y\|_{\psi_1}. \tag{B.4}$$

Therefore, the space of subexponential random variables forms a normed vector space. There are many equivalent norms for subexponential and subgaussian random variables, see more in §5.2.3 and §5.2.4 in Vershynin (2010).

We have the following lemma that relates subgaussian and subexponential norms.

**Lemma 6.** For subexponential random variables $X$ and $Y$ (*not* necessarily independent) we have

$$\|XY\|_{\psi_1} \leq \|X\|_{\psi_2}^2 + \|Y\|_{\psi_2}^2. \tag{B.5}$$

*Proof.* (i) We first prove

$$\|X^2\|_{\psi_1} \leq 2\|X\|_{\psi_2}^2. \tag{B.6}$$

From the definition of subgaussian and subexponential norm in Definition 3, note for each $p_* \geq 1$

$$p_*^{-1} \left(\mathbb{E}|X|^{2p_*}\right)^{1/p_*} \leq 2 \cdot (2p_*)^{-1} \left(\mathbb{E}|X^2|^{p_*}\right)^{1/p_*} \leq 2\|X\|_{\psi_2}^2.$$

Taking sup over $p_* \geq 1$ we complete the proof of (B.6).

(ii) Note first $2|XY| \leq X^2 + Y^2$, from Definition 3 we have

$$2\|XY\|_{\psi_1} = \|2|XY|\|_{\psi_1} \leq \|X^2 + Y^2\|_{\psi_1}$$

To conclude (B.5), note by triangle inequality (B.4) and (B.6)

$$\|XY\|_{\psi_1} \leq \frac{1}{2}\|X^2 + Y^2\|_{\psi_1} \leq \frac{1}{2}\left(\|X^2\|_{\psi_1} + \|Y^2\|_{\psi_1}\right)$$
$$\leq \frac{1}{2}\left(2\|X\|_{\psi_2}^2 + 2\|Y\|_{\psi_2}^2\right) \leq \|X\|_{\psi_2}^2 + \|Y\|_{\psi_2}^2.$$

and we are done. □

**Lemma 7** (Moment generating function of subexponential random variables). Let $X$ be a subexponential random variable with $\mathbb{E}X = 0$, so $\|X\|_{\psi_1} < \infty$. Then for $t$ such that $|t| \leq (C_{7,1}\|X\|_{\psi_1})^{-1}$

$$\mathbb{E}\exp(tX) \leq \exp\left(C_{7,2}t^2\|X\|_{\psi_1}^2\right).$$

In particular, there is a constant $K = C_{7,1}\|X\|_{\psi_1}$ such that

$$\mathbb{E}\exp(tX) \leq \exp\left(\frac{t^2K^2}{2}\right).$$

*Proof.* The case where $\|X\|_{\psi_1} = 0$ is trivial, since $\mathbb{E}X = 0$ implies $X = 0, a.s.$, and both sides are equal to 1. Assume without loss of generality that $\|X\|_{\psi_1} = 1$; for the general case the result is obtained by writing $t\|X\|_{\psi_1}$ in the place of $t$ and $X/\|X\|_{\psi_1}$ in the place of $X$. Therefore

$$\mathbb{E}|X|^p \leq p^p\|X\|_{\psi_1} \leq p^p.$$

Taylor expansion and $\mathbb{E}X = 0$ together indicates

$$\mathbb{E}\exp(tX) = \sum_{p=0}^{\infty}\frac{\mathbb{E}X^p}{p!}t^p \leq 1 + \sum_{p=2}^{\infty}\frac{\mathbb{E}|X|^p}{p!}|t|^p \leq 1 + \sum_{p=2}^{\infty}\frac{p^p}{p!}|t|^p. \tag{B.7}$$

An application of Stirling's formula allows that for each $p \geq 1$

$$\frac{p^p}{p!} \leq \frac{1}{\sqrt{2\pi p}}\cdot e^p \leq e^p.$$

If $|t| \leq 1/(2e)$ then the right hand side of (B.7) has the following upper bound

$$\mathbb{E}\exp(tX) \leq 1 + \sum_{p=2}^{\infty}(e|t|)^p \leq 1 + 2e^2t^2 \leq \exp(2e^2t^2).$$

Taking $C_{7,1} = 2e$ and $C_{7,2} = 2e^2$ completes the proof.

□

## B.3 A Bernoulli-type of concentration inequality for subexponential martingale difference sequences

In this subsection we introduce Lemma 8, a Bernoulli-type of concentration inequality for martingale difference sequences whose tails are *subexponential*. Since we cannot find the proof in any references we give all the proof details, which is of independent interest.

**Lemma 8** (Bernstein-type concentration inequality). Let $\xi_1, \ldots, \xi_n$ be a martingale difference sequence with regards to the filtration $\mathcal{F}_n$, and furthermore *uniformly subexponential* in the sense that there is a $K > 0$ satisfying for each $1 \leq i \leq n$ and $|t| \leq K^{-1}$

$$\mathbb{E}\left[\exp(t\xi_i) \mid \mathcal{F}_{i-1}\right] \leq \exp\left(\frac{t^2K^2}{2}\right), \quad a.s.$$

Then for each $\mathbf{a} = (a_1, \ldots, a_n) \in \mathbb{R}^n$, $z > 0$ and $n \geq 1$

$$\mathbb{P}\left(\left|\sum_{i=1}^{n}a_i\xi_i\right| \geq z\right) \leq 2\exp\left(-C_8\min\left(\frac{z^2}{K^2\|\mathbf{a}\|_2^2}, \frac{z}{K\|\mathbf{a}\|_\infty}\right)\right).$$

*Proof.* Assume without loss of generality that $K = 1$; for general $K$ the result is obtained by writing $X_i/K$ in the place of $X_i$ and $t/K$ in the place of $t$. Since for $0 < t \leq \|\mathbf{a}\|_\infty^{-1}$ we have $t|a_i| \leq 1$ for $1 \leq i \leq n$ and hence

$$\mathbb{E}\left[\exp\left(t\sum_{i=1}^{n} a_i\xi_i\right) \Big| \mathcal{F}_{n-1}\right] \leq \exp\left(t\sum_{i=1}^{n-1} a_i\xi_i\right) \mathbb{E}\left[\exp(ta_n\xi_n) \mid \mathcal{F}_{n-1}\right]$$

$$\leq \exp\left(t\sum_{i=1}^{n-1} a_i\xi_i\right) \exp\left(C_{7,1}^{-2}C_{7,2}t^2 a_n^2\right).$$

Taking expectation and using the law of iterated expectation

$$\mathbb{E}\left[\exp\left(t\sum_{i=1}^{n} a_i\xi_i\right)\right] \leq \mathbb{E}\left[\exp\left(t\sum_{i=1}^{n-1} a_i\xi_i\right)\right] \exp\left(C_{7,1}^{-2}C_{7,2}t^2 a_n^2\right).$$

Applying the above formula iteratively we obtain

$$\mathbb{E}\left[\exp\left(t\sum_{i=1}^{n} a_i\xi_i\right)\right] \leq \exp\left(C_{7,1}^{-2}C_{7,2}t^2 \sum_{i=1}^{n} a_i^2\right) = \exp\left(C_{7,1}^{-2}C_{7,2}t^2 \|\mathbf{a}\|_2^2\right).$$

Now apply Markov's inequality and use the above display we conclude for each $z > 0$

$$\mathbb{P}\left(\sum_{i=1}^{n} a_i\xi_i \geq z\right) = \mathbb{P}\left(\exp\left(t\sum_{i=1}^{n} a_i\xi_i\right) \geq \exp(tz)\right)$$

$$\leq \exp(-tz)\, \mathbb{E}\left[\exp\left(t\sum_{i=1}^{n} a_i\xi_i\right)\right] \leq \exp\left(-tz + C_{7,1}^{-2}C_{7,2}t^2 \|\mathbf{a}\|_2^2\right).$$

We optimize the right hand side and plug in $t^* = \min\left((2C_{7,1}^{-2}C_{7,2})^{-1}\|\mathbf{a}\|_2^{-2} z, \|\mathbf{a}\|_\infty^{-1}\right)$ in the place of $t$ to obtain

$$\mathbb{P}\left(\sum_{i=1}^{n} a_i\xi_i \geq z\right) \leq \exp\left(-tz + C_{7,1}^{-2}C_{7,2}t^2 \|\mathbf{a}\|_2^2\right)$$

$$\leq \max\left(\exp\left(-\frac{z^2}{4C_{7,1}^{-2}C_{7,2}\|\mathbf{a}\|_2^2}\right), \exp\left(-\frac{z}{2\|\mathbf{a}\|_\infty}\right)\right)$$

$$\leq \exp\left(-C_8 \min\left(\frac{z^2}{\|\mathbf{a}\|_2^2}, \frac{z}{\|\mathbf{a}\|_\infty}\right)\right),$$

where $C_8 \equiv \left(4C_{7,1}^{-2}C_{7,2}\right)^{-1} \wedge (1/2)$. Repeat the above argument for $-\xi_i$'s in the place of $\xi_i$ we obtain the same bound for $\mathbb{P}\left(-\sum_{i=1}^{n} a_i\xi_i \geq z\right)$, which completes the proof altogether. □

## C  Proof of Auxiliary Propositions

This section provides the proofs of auxiliary propositions. For brevity in this section, we let $\mathcal{F}_n = \sigma(\mathbf{Y}^{(i)} : i = 0, 1, \ldots, n)$ be the filtration we are focused on, i.e., the smallest $\sigma$-algebra with regards to which the initial $n$ samples are measurable.

## C.1 Proof of Proposition 1

We first conduct an analysis on the increment of $U_k^{(n)}$ defined in (3.8). Recall that we initialize our algorithm from $\mathbf{v}^{(0)} \in \mathcal{S}_2$. Recall (3.6), and we conclude the following lemma.

**Lemma 9.** Assume all conditions in Theorem 1 along with the warm start condition $\mathbf{v}^{(0)} \in \mathcal{S}_2$. For each $k = 2, \ldots, d$ and any $n \geq 1$, there exists a random variable $R_{k,n}$ satisfying that on the event $(\mathcal{N}_M > n, \mathcal{N}_w \geq n)$ we have

$$|R_{k,n}| \leq C_9 M^2 \beta^2, \quad a.s., \tag{C.1}$$

where $\mathcal{N}_M$ was defined earlier in (3.6), such that the increment of $U_k^{(n)}$ has

$$U_k^{(n)} - U_k^{(n-1)} = \beta \cdot \frac{\mathbf{v}^\top \mathbf{Y}}{v_1} \cdot \left(Y_k - \frac{v_k}{v_1} Y_1\right) + R_{k,n}. \tag{C.2}$$

*Proof.* To derive (C.2), note we have

$$U_k^{(n)} - U_k^{(n-1)} = \frac{v_k + \beta(\mathbf{v}^\top \mathbf{Y})Y_k}{v_1 + \beta(\mathbf{v}^\top \mathbf{Y})Y_1} - \frac{v_k}{v_1} = \beta\left(1 + \beta\frac{\mathbf{v}^\top \mathbf{Y}}{v_1}Y_1\right)^{-1} \frac{\mathbf{v}^\top \mathbf{Y}}{v_1}\left(Y_k - \frac{v_k}{v_1}Y_1\right)$$

$$= \beta\frac{\mathbf{v}^\top \mathbf{Y}}{v_1}\left(Y_k - \frac{v_k}{v_1}Y_1\right) + R_{k,n},$$

where

$$R_{k,n} = \beta\left[\left(1 + \beta\frac{\mathbf{v}^\top \mathbf{Y}}{v_1}Y_1\right)^{-1} - 1\right] \cdot \frac{\mathbf{v}^\top \mathbf{Y}}{v_1} \cdot \left(Y_k - \frac{v_k}{v_1}Y_1\right).$$

We note for $\mathbf{v} \in \mathcal{S}_3$, on $(\mathcal{N}_M > n)$

$$\left|\beta\frac{\mathbf{v}^\top \mathbf{Y}}{v_1}Y_1\right| \leq 3M\beta \leq 1/2.$$

Since $\left|(1+x)^{-1} - 1\right| = |x||1+x|^{-1} \leq 2|x|$ for $|x| \leq \frac{1}{2}$, we have

$$|R_{k,n}| \leq 2\beta^2 \left(\frac{\mathbf{v}^\top \mathbf{Y}}{v_1}\right)^2 |Y_1| \left|Y_k - \frac{v_k}{v_1}Y_1\right| \leq CM^2\beta^2.$$

This concludes (C.2). □

Our strategy to study the $U_k^{(n)}$ is as follows: first by subtracting off the drift terms one can extract out a martingale difference sequence, and then use the martingale concentration inequality to bound the martingale difference sequence using the subgaussian bounds.

**Lemma 10.** On $(\mathcal{N}_w \geq n) \in \mathcal{F}_{n-1}$ we have

$$\mathbb{E}\left[\beta\frac{\mathbf{v}^\top \mathbf{Y}}{v_1}\left(Y_k - \frac{v_k}{v_1}Y_1\right) \bigg| \mathcal{F}_{n-1}\right] = -\beta(\lambda_1 - \lambda_k)U_k^{(n-1)}. \tag{C.3}$$

*Proof.* Note

$$\mathbb{E}\left[\beta\frac{\mathbf{v}^\top \mathbf{Y}}{v_1}\left(Y_k - \frac{v_k}{v_1}Y_1\right)\bigg|\mathcal{F}_{n-1}\right] = \beta\frac{1}{v_1}\left(\lambda_k v_k - \frac{v_k}{v_1}\lambda_1 v_1\right) = -\beta\left(\lambda_1 - \lambda_k\right)\frac{v_k}{v_1}.$$

Since $U_k^{(n-1)} = v_k^{(n-1)}/v_1^{(n-1)}$ we conclude (C.3).

□

With Lemma 10 and (C.2) in Lemma 9 at hand, let us define for each $n = 1, 2, \ldots$

$$e_{k,n} = \frac{\mathbf{v}^\top \mathbf{Y}}{v_1}\left(Y_k - \frac{v_k}{v_1}Y_1\right) + \beta(\lambda_1 - \lambda_k)U_k^{(n-1)}. \tag{C.4}$$

From (C.4) we can express the $U_k^{(n)}$ update rule as

$$U_k^{(n)} = e_{k,n} + (1 - \beta(\lambda_1 - \lambda_k))U_k^{(n-1)} + R_{k,n}. \tag{C.5}$$

In parallel, let

$$\bar{R}_{k,n} = R_{k,n}\mathbf{1}_{(\mathcal{N}_M > n, \mathcal{N}_w \geq n)}, \tag{C.6}$$

let

$$\bar{e}_{k,n} = e_{k,n}\mathbf{1}_{(\mathcal{N}_w \geq n)}, \tag{C.7}$$

let $\bar{U}_k^{(0)} = U_k^{(0)}$, and define the coupled process iteratively

$$\bar{U}_k^{(n)} = \bar{e}_{k,n} + (1 - \beta(\lambda_1 - \lambda_k))\bar{U}_k^{(n-1)} + \bar{R}_{k,n}. \tag{C.8}$$

We conclude the following lemma that characterize the coupling relations.

**Lemma 11.** *For each $n \geq 0$ and $k = 2, \ldots, d$ we have $\bar{U}_k^{(n)} = U_k^{(n)}$ on the event $(\mathcal{N}_M > n, \mathcal{N}_w \geq n)$.*

*Proof.* For $n = 0$ the lemma holds from definition. In general if it holds for $n - 1$ then since by applying (C.5) and (C.8):

From the definition of $\bar{R}_{k,i}$ and $e_{k,i}$ in (C.6) and (C.7) we have on $(\mathcal{N}_M > n, \mathcal{N}_w \geq n)$ that $\bar{R}_{k,i} = R_{k,i}$ and $\bar{e}_{k,i} = e_{k,i}$ for all $i \leq n$ and hence conclude the lemma.

□

The coupling relation in Lemma 11 allows us to analyze the $\bar{U}_k^{(n)}$ iteration, which enjoys desirable tail properties and concentration inequalities. Let

$$E_{k,n-1} = \mathbb{E}\left[\bar{R}_{k,n} \mid \mathcal{F}_{n-1}\right],$$

and

$$f_{k,n} = \bar{e}_{k,n} + \bar{R}_{k,n} - E_{k,n-1}, \tag{C.9}$$

and from (C.1) we have $|E_{k,n-1}| \leq |\bar{R}_{k,n}| \leq C_9 M^2 \beta^2$, a.s. From (C.8)

$$\bar{U}_k^{(n)} = f_{k,n} + (1 - \beta(\lambda_1 - \lambda_k))\bar{U}_k^{(n-1)} + E_{k,n-1}. \tag{C.10}$$

Iteratively applying (C.10) we have

$$\bar{U}_k^{(n)} = (1 - \beta(\lambda_1 - \lambda_k))^n U_k^{(0)} + \sum_{i=1}^n (1 - \beta(\lambda_1 - \lambda_k))^{n-i} f_{k,i} \\ + \sum_{i=1}^n (1 - \beta(\lambda_1 - \lambda_k))^{n-i} E_{k,i-1}. \quad (C.11)$$

Since the event $(\mathcal{N}_w < i)$ and its complement $(\mathcal{N}_w \geq i)$ are both measurable to $\mathcal{F}_{i-1}$, we have $\mathbb{E}[\bar{e}_{k,i} \mid \mathcal{F}_{i-1}] = \mathbb{E}[e_{k,i} \mid \mathcal{F}_{i-1}]1_{\mathcal{N}_w \geq i} = 0$ so it also forms a martingale sequence. We have for any fixed $n \geq 1$ and $k = 2, \ldots, d$ that $(f_{k,i} : 1 \leq i \leq n)$ defined in (C.9) forms a martingale difference sequence. Hence one can turn to analyze the martingale difference sequence $f_{k,n}$, which enjoys *conditional subexponential tails* and has the following concentration inequality.

**Lemma 12.** We have for all $y$ such that

$$0 < y \leq C_{12,1} \cdot \widehat{\beta}^{-0.5+\varepsilon} \leq C_{12,1}\lambda_1^{-1} \cdot \widehat{\beta}^{0.5+\varepsilon}\beta^{-1}, \quad (C.12)$$

the following concentration result holds

$$\mathbb{P}\left(\left|\sum_{i=1}^n (1 - \beta(\lambda_1 - \lambda_k))^{n-i} f_{k,i}\right| \geq y\widehat{\beta}^{0.5-\varepsilon}\right) \leq 2\exp\left(-C_{12,2}y^2\widehat{\beta}^{-2\varepsilon}\right). \quad (C.13)$$

*Proof.* (i) We first bound the subexponential norms standing at $\mathcal{F}_{n-1}$. Using Lemma 6 we conclude from Assumption 1

$$\|(\mathbf{v}^\top \boldsymbol{Y})Y_k\|_{\psi_1} \leq \|\mathbf{v}^\top \boldsymbol{Y}\|_{\psi_2}^2 + \|Y_k\|_{\psi_2}^2 \leq \lambda_1 + \lambda_k \leq 2\lambda_1.$$

From the definition of $\mathcal{S}_3$ in (3.4), $|(v_k/v_1)1_{\mathcal{N}_w \geq n}| \leq 3$ and

$$\|(v_k/v_1)(\mathbf{v}^\top \boldsymbol{Y})Y_1 1_{\mathcal{N}_w \geq n}\|_{\psi_1} \leq 3\left(\|\mathbf{v}^\top \boldsymbol{Y}\|_{\psi_2}^2 + \|Y_1\|_{\psi_2}^2\right) \leq 6\lambda_1,$$

where we used $\|Y_k\|_{\psi_2} \leq \lambda_k^{1/2}$ for each $k$ and $\|\mathbf{v}^\top \boldsymbol{Y}\|_{\psi_2} \leq \lambda_1^{1/2}$. Hence

$$\left\|\frac{\beta}{v_1}\left((\mathbf{v}^\top \boldsymbol{Y})Y_k - \frac{v_k}{v_1}(\mathbf{v}^\top \boldsymbol{Y})Y_1\right)1_{\mathcal{N}_w \geq n}\right\|_{\psi_1} \leq 3\beta(8\lambda_1) = 24\lambda_1\beta.$$

Subtracting off its conditional expectation we conclude that $\|\bar{e}_{k,n}\|_{\psi_1} \leq 2(24\lambda_1) = 48\lambda_1$. For the residual term we have from $M^2\beta^2 \leq \lambda_1\beta$ that

$$\left\|\bar{R}_{k,n} - E_{k,n-1}\right\|_{\psi_1} \leq 2C_9 M^2\beta^2 \leq 2C_9\lambda_1\beta.$$

From (C.9) we have

$$\|f_{k,n}\|_{\psi_1} \leq \|\bar{e}_{k,n}\|_{\psi_1} + \left\|\bar{R}_{k,n} - E_{k,n-1}\right\|_{\psi_1} \leq (48 + 2C_9)\lambda_1\beta.$$

Letting $C_* = 48 + 2C_9$ and applying Lemma 7 allows us to conclude for $K = C_{7,1}C_*\lambda_1\beta$ and $|t| \leq K^{-1}$

$$\mathbb{E}\left[\exp(tf_{k,n}) \mid \mathcal{F}_{n-1}\right] \leq \exp\left(\frac{t^2 K^2}{2}\right). \quad (C.14)$$

(ii) We now come to prove the concentration inequality. Let for $i = 1, \ldots, n$ $a_i = (1 - \beta(\lambda_1 - \lambda_k))^{n-i}$ as the $i$th coordinate of $\mathbf{a}$, and also $f_{k,i}$ in the place of $\xi_i$, then the sum

$$\|\mathbf{a}\|_2^2 = \sum_{i=1}^n (1 - \beta(\lambda_1 - \lambda_k))^{2(n-i)} \leq \frac{1}{2\beta(\lambda_1 - \lambda_k) - \beta^2(\lambda_1 - \lambda_k)^2} \leq \frac{1}{\beta(\lambda_1 - \lambda_k)}, \quad \text{(C.15)}$$

where we used $\beta(\lambda_1 - \lambda_k) \leq 1$. Applying Lemma 8 with

$$\|\mathbf{a}\|_\infty = 1, \qquad \|\mathbf{a}\|_2 = (\lambda_1 - \lambda_k)^{-0.5}\beta^{-0.5}, \qquad K = C_*\lambda_1\beta, \qquad z = y\widehat{\beta}^{0.5-\varepsilon},$$

we have from $\widehat{\beta} \geq \lambda_1^2(\lambda_1 - \lambda_k)^{-1}\beta$

$$\frac{z^2}{K^2\|\mathbf{a}\|_2^2} = \frac{y^2\widehat{\beta}^{1-2\varepsilon}}{C_*^2\lambda_1^2\beta^2 \cdot (\lambda_1 - \lambda_k)^{-1}\beta^{-1}} = C_*^{-2}y^2\widehat{\beta}^{1-2\varepsilon} \cdot (\lambda_1 - \lambda_k)\lambda_1^{-2} \cdot \beta^{-1} \geq C_*^{-2}y^2\widehat{\beta}^{-2\varepsilon}, \quad \text{(C.16)}$$

and

$$\frac{z}{K\|\mathbf{a}\|_\infty} = \frac{y\widehat{\beta}^{0.5-\varepsilon}}{C_*\lambda_1\beta \cdot 1} = C_*^{-1}y \cdot \lambda_1^{-1} \cdot \widehat{\beta}^{0.5-\varepsilon}\beta^{-1}, \quad \text{(C.17)}$$

so when (C.12) holds, the minimum of (C.16) and (C.17) is

$$\min\left(\frac{z^2}{K^2\|\mathbf{a}\|_2^2}, \frac{z}{K\|\mathbf{a}\|_\infty}\right) \geq \min\left(C_*^{-2}y^2\widehat{\beta}^{-2\varepsilon}, C_*^{-1}y \cdot \lambda_1^{-1} \cdot \widehat{\beta}^{0.5-\varepsilon}\beta^{-1}\right) = C_*^{-2}y^2\widehat{\beta}^{-2\varepsilon}.$$

Hence using (C.15) and the concentration inequality as in Lemma 8 we prove using (C.14) that for all $y$ satisfying (C.12) and $n \geq 1$

$$\mathbb{P}\left(\left|\sum_{i=1}^n (1 - \beta(\lambda_1 - \lambda_k))^{n-i} f_{k,i}\right| \geq y\widehat{\beta}^{0.5-\varepsilon}\right) \leq 2\exp\left(-C_{8,3}C_*^{-2} \cdot y^2\widehat{\beta}^{-2\varepsilon}\right). \quad \text{(C.18)}$$

Setting in (C.18) $C_{12,1} = C_*$ and $C_{12,2} = C_{8,3}C_*^{-2}$ proves (C.13) and hence the lemma. □

We are now ready to prove the Proposition 1.

*Proof of Proposition 1.* To bound the sum in the first line of (C.11), by noticing $\beta(\lambda_1 - \lambda_k) < 1$ and the fact $|E_{k,i-1}| \leq C_9 M^2\beta^2$ we have

$$\left|\sum_{i=1}^n (1 - \beta(\lambda_1 - \lambda_k))^{n-i} E_{k,i-1}\right| \leq C_9 M^2\beta^2 \sum_{i=1}^n (1 - \beta(\lambda_1 - \lambda_k))^{n-i}$$

$$\leq \frac{C_9 M^2\beta}{\lambda_1 - \lambda_k} \leq C_9\widehat{\beta}. \quad \text{(C.19)}$$

Moreover for the sum in the second line of (C.11), used $\beta(\lambda_1 - \lambda_k) \leq \beta(\lambda_1 - \lambda_2) \leq 1$ which is guaranteed if $\widehat{\beta} \leq 1$ which is further guaranteed if $d\widehat{\beta}^{0.5} \leq 1$ (C.13) implies when $C_9 \leq C_{12,1} \cdot \widehat{\beta}^{-0.5+\varepsilon} \leq C_{12,1}\lambda_1^{-1} \cdot \widehat{\beta}^{0.5+\varepsilon}\beta^{-1}$ i.e. $\widehat{\beta}^{1-2\varepsilon} \leq (C_{12,1} \cdot C_9^{-1})^2 \wedge 1$,

$$\mathbb{P}\left(\left|\sum_{i=1}^n (1 - \beta(\lambda_1 - \lambda_k))^{n-i} f_{k,i}\right| \geq C_9\widehat{\beta}^{0.5-\varepsilon}\right) \leq 2\exp\left(-C_9^2 C_{12,2}\widehat{\beta}^{-2\varepsilon}\right).$$

Combining (C.11), (C.13) and (C.19) one has when $\widehat{\beta} \leq 1$

$$\mathbb{P}\left(\sup_{n \leq N_t^* \wedge \mathcal{N}_w} \left|\bar{U}_k^{(n)} - U_k^{(0)}(1 - \beta(\lambda_1 - \lambda_k))^n\right| \geq 2C_9 \widehat{\beta}^{0.5-\varepsilon}\right) \leq 2N_t^* \exp\left(-C_9^2 C_{12,2} \widehat{\beta}^{-2\varepsilon}\right).$$

Setting $C'_{1,P} = C_9^2 C_{12,2}$ and $C'_{1,P} = 2C_9$ and noting $(\mathcal{N}_M > N) \subseteq \left(\bar{U}_k^{(n)} = U_k^{(n)}, n \leq N_{\beta,t}^* \wedge \mathcal{N}_w\right)$ from Lemma 11 completes the proof of Proposition 1. □

## C.2 Proof of Proposition 2

Let the stopping time

$$\mathcal{N}_v = \inf\left\{n \geq N_{0.5-\varepsilon}^* : \sup_{k:2 \leq k \leq d} \left|U_k^{(n)}\right| > 2\widehat{\beta}^{0.5-\varepsilon}\right\}, \tag{C.20}$$

where by convention, $\inf \varnothing = \infty$ and hence $\mathcal{N}_v$ is always no less than $N_{0.5-\varepsilon}^*$. We first give a probability estimate of $(\mathcal{N}_v \leq N_{\beta,t}^*, \mathcal{N}_M > N_{\beta,t}^*)^c$ which occurs with high probability under our assumptions.

**Lemma 13.** We have

$$\mathbb{P}\left(\mathcal{N}_v \leq N_{\beta,t}^*, \mathcal{N}_M > N_{\beta,t}^*\right) \leq 2(d-1)N_{\beta,t}^* \exp\left(-C_{13}\widehat{\beta}^{-2\varepsilon}\right). \tag{C.21}$$

*Proof.* For each $k = 2, \ldots, d$ we define the events

$$\mathcal{J}_{1,k} \equiv \left(\sup_{n \leq N_{\beta,t}^* \wedge \mathcal{N}_w} \left|U_k^{(n)} - U_k^{(0)}(1 - \beta(\lambda_1 - \lambda_k))^n\right| > \widehat{\beta}^{0.5-\varepsilon}\right) \cap (\mathcal{N}_M > N_{\beta,t}^*),$$

$$\mathcal{J}_1 \equiv \bigcup_{k=2}^d \mathcal{J}_{1,k}. \tag{C.22}$$

Taking $t = 2$ in Proposition 1 implies that for $k = 2, \ldots, d$, $\mathbb{P}(\mathcal{J}_{1,k}) \leq 2N_{\beta,t}^* \exp\left(-C_{1,P}\widehat{\beta}^{-2\varepsilon}\right)$, we conclude

$$\mathbb{P}(\mathcal{J}_1) \leq \sum_{k=2}^d \mathbb{P}\left(\mathcal{J}_{1,k}^c\right) \leq 2(d-1)N_{\beta,t}^* \exp\left(-C_{1,P}\widehat{\beta}^{-2\varepsilon}\right),$$

Therefore in order to prove the lemma one only needs to prove $(\mathcal{N}_v \leq N_{\beta,t}^*, \mathcal{N}_M > N_{\beta,t}^*) \subseteq \mathcal{J}_1$ which is equivalent to

$$\mathcal{J}_1^c \cap (\mathcal{N}_M > N_{\beta,t}^*) \subseteq (\mathcal{N}_v > N_{\beta,t}^*). \tag{C.23}$$

The rest of this proof devotes to prove (C.23).

Since on the event $\mathcal{J}_1^c \cap (\mathcal{N}_M > N_{\beta,t}^*, \mathcal{N}_w < N_{\beta,t}^*)$, $N_{\beta,t}^* \wedge \mathcal{N}_w = \mathcal{N}_w$ so

$$\sum_{k=2}^d \left|U_k^{(\mathcal{N}_w)} - U_k^{(0)}(1 - \beta(\lambda_1 - \lambda_k))^{\mathcal{N}_w}\right|^2 \leq (d-1)\widehat{\beta}^{1-2\varepsilon}.$$

By triangle inequality for Euclidean norm and $d\widehat{\beta}^{1-2\varepsilon} \leq 1$

$$\left(\sum_{k=2}^{d} |U_k^{(\mathcal{N}_w)}|^2\right)^{1/2} \leq \left(\sum_{k=2}^{d} |U_k^{(0)}|^2\right)^{1/2} + \left((d-1)\widehat{\beta}^{1-2\varepsilon}\right)^{1/2} \leq 1 + \left(d\widehat{\beta}^{1-2\varepsilon}\right)^{1/2} \leq 2.$$

Nevertheless the above display contradicts with the definition of $\mathcal{N}_w$ in (3.5) and $\mathcal{S}_3$ in (3.4) where

$$\sum_{k=2}^{d} |U_k^{(\mathcal{N}_w)}|^2 = \left(v_1^{(\mathcal{N}_w)}\right)^{-2} - 1 \geq 8.$$

This indicates $(\mathcal{N}_M > N_{\beta,t}^*, \mathcal{N}_w < N_{\beta,t}^*) \cap \mathcal{J}_1^c = \emptyset$, i.e. $\mathcal{J}_1^c \cap (\mathcal{N}_M > N_{\beta,t}^*) \subseteq (\mathcal{N}_w \geq N_{\beta,t}^*)$.

From the previous analysis, by definition in (C.22) we know on the event $\mathcal{J}_1^c \cap (\mathcal{N}_M > N_{\beta,t}^*)$, for each $k = 2, \ldots, d$

$$\sup_{n \leq N_{\beta,t}^*} \left|U_k^{(n)} - U_k^{(0)}(1 - \beta(\lambda_1 - \lambda_k))^n\right| \leq \widehat{\beta}^{0.5-\varepsilon},$$

and hence for $n \in [N_{\beta,0.5-\varepsilon}^*, N_{\beta,t}^*]$ and each $k = 2, \ldots, d$, $(1 - \beta(\lambda_1 - \lambda_k))^{N_{\beta,0.5-\varepsilon}^*} \leq \widehat{\beta}^{0.5-\varepsilon}$, indicating

$$\left|U_k^{(n)}\right| \leq \widehat{\beta}^{0.5-\varepsilon} + \left|U_k^{(0)}\right|(1 - \beta(\lambda_1 - \lambda_k))^n \leq \widehat{\beta}^{0.5-\varepsilon} + \widehat{\beta}^{0.5-\varepsilon} = 2\widehat{\beta}^{0.5-\varepsilon}.$$

From the definition in (C.20) we know $(\mathcal{N}_v > N_{\beta,t}^*)$ occurs, and hence conclude the right equation of (C.23), completing the proof.

$\square$

The intuition is that, from (C.21) we have $(\mathcal{N}_v > n)$ occurs with high probability for $n \leq \mathcal{N}_{\beta,t}^*$. In the following all the probability measure are conditional on $\mathcal{F}_{n-1}$. Starting from $\mathcal{N}_{\beta,0.5-\varepsilon}^*$ one can bound for each $k = 2, \ldots, d$ $U_k^{(n)}$ within a very small neighborhood of 0. Let for each $n = 1, 2, \ldots$

$$h_{k,n} = f_{k,n}\mathbf{1}_{(\mathcal{N}_v \geq n)}, \tag{C.24}$$

where $f_{k,n}$ is defined in (C.9). We conclude the following facts:

(i) $h_{k,n} = f_{k,n}$ for $n \leq \mathcal{N}_{\beta,0.5-\varepsilon}^*$. This is because $(\mathcal{N}_v \geq N_{0.5-\varepsilon}^*)$ always occurs;

(ii) $h_{k,n}$ forms a martingale difference sequence with regards to filtration $\mathcal{F}_n$, due to $(\mathcal{N}_v \geq n) \in \mathcal{F}_{n-1}$.

We first obtain a tight estimate for the second moment of error term $h_{k,n}$.

**Lemma 14.** For $n \leq N_{\beta,0.5-\varepsilon}^*$ we have

$$\mathbb{E}\left[h_{k,n}^2 \mid \mathcal{F}_{n-1}\right] \leq 16\beta^2 v_1^{-2}\left(\lambda_1\lambda_k + \left(\frac{v_k}{v_1}\right)^2 \lambda_1^2\right) \leq C_{14,1}\beta^2\lambda_1^2. \tag{C.25}$$

For $n \in \left[N_{\beta,0.5-\varepsilon}^* + 1, N_{\beta,t}^*\right]$ a sharper upper bound is available

$$\mathbb{E}\left[h_{k,n}^2 \mid \mathcal{F}_{n-1}\right] \leq C_{14,2}\beta^2\left(\lambda_1\lambda_k + 4\widehat{\beta}^{1-2\varepsilon}\lambda_1^2\right). \tag{C.26}$$

*Proof.* (C.2) in Lemma 9 implies that

$$\mathbb{E}[h_{k,n}^2] = \mathbf{1}_{(\mathcal{N}_v \geq n)} \mathbb{E}[f_{k,n}^2] = \mathbf{1}_{(\mathcal{N}_v \geq n)} \mathrm{var}\left(\beta \frac{\mathbf{v}^\top \mathbf{Y}}{v_1}\left(Y_k - \frac{v_k}{v_1} Y_1\right)\right).$$

We have by Cauchy-Schwarz inequality that for each $k = 1, \ldots, d$:

$$\mathrm{var}\left[\beta \frac{\mathbf{v}^\top \mathbf{Y}}{v_1}\left(Y_k - \frac{v_k}{v_1} Y_1\right)\right] \leq \frac{\beta^2}{v_1^2} \mathbb{E}\left[(\mathbf{v}^\top \mathbf{Y})^2 \left(Y_k - \frac{v_k}{v_1} Y_1\right)^2\right]$$

$$\leq \frac{\beta^2}{v_1^2} \left[\mathbb{E}(\mathbf{v}^\top \mathbf{Y})^4\right]^{1/2} \left[\mathbb{E}\left(Y_k - \frac{v_k}{v_1} Y_1\right)^4\right]^{1/2}.$$

Here $\mathbf{Y} = \mathbf{\Lambda}^{1/2} \mathbf{Z}$, and for each $\mathbf{v} \in \mathcal{S}^{d-1}$ we have from (B.3) in Lemma 5 that

$$\mathbb{E}(\mathbf{v}^\top \mathbf{Y})^4 = \mathbb{E}(\mathbf{v}^\top \mathbf{\Lambda}^{1/2} \mathbf{Z})^4 \leq 16 \|\mathbf{\Lambda}^{1/2} \mathbf{v}\|^4 \leq 16 \lambda_1^2,$$

and from (C.27)

$$\mathbb{E}\left(Y_k - \frac{v_k}{v_1} Y_1\right)^4 = \mathbb{E}\left(\lambda_k^{1/2} Z_k - \frac{v_k}{v_1} \lambda_1^{1/2} Z_1\right)^4 \leq 16 \left(\lambda_k + \left(\frac{v_k}{v_1}\right)^2 \lambda_1\right)^2.$$

Combining the last three displays concludes (C.25). For (C.26) we have on $(\mathcal{N}_v \geq n)$

$$\left(\frac{v_k}{v_1}\right)^2 \leq 4\widehat{\beta}^{1-2\varepsilon}, \qquad v_1^{-2} \leq 9, \tag{C.27}$$

so when $n \geq N^*_{\beta, 0.5-\varepsilon}$

$$\mathbf{1}_{(\mathcal{N}_v \geq n)} \cdot \mathbb{E}\left(Y_k - \frac{v_k}{v_1} Y_1\right)^4 \leq 16 \left(\lambda_k + 4\widehat{\beta}^{1-2\varepsilon} \lambda_1\right)^2.$$

From (C.27) we have

$$\mathbf{1}_{(\mathcal{N}_v \geq n)} \cdot \mathbb{E}\left[\beta \frac{\mathbf{v}^\top \mathbf{Y}}{v_1}\left(Y_k - \frac{v_k}{v_1} Y_1\right)\right]^2 \leq \mathbf{1}_{(\mathcal{N}_v \geq n)} \cdot 16\beta^2 v_1^{-2} \left(\lambda_1 \lambda_k + 4\widehat{\beta}^{1-2\varepsilon} \lambda_1^2\right)$$

$$\leq 16\beta^2 \cdot 9 \cdot \left(\lambda_1 \lambda_k + 4\widehat{\beta}^{1-2\varepsilon} \lambda_1^2\right).$$

Letting $C_{14,2} = 16 \cdot (9)$ completes the proof. $\square$

Recall that the truncated noise terms $h_{k,i}$ were defined in (C.24), and we let

$$\widehat{W}_k^{(n)} = \sum_{i=1}^{n} (1 - \beta(\lambda_1 - \lambda_k))^{n-i} h_{k,i}. \tag{C.28}$$

From the viewpoint of (C.11), we couple $U_k^{(n)}, n > N^*_{\beta, 0.5-\varepsilon}$ with a sequence in parallel that uses $h_{k,i}$'s

$$W_k^{(n)} = (1 - \beta(\lambda_1 - \lambda_k))^n U_k^{(0)} + \widehat{W}_k^{(n)}. \tag{C.29}$$

We prove the following Lemma 15.

**Lemma 15.** We have for each $n \geq N^*_{\beta,1-2\varepsilon}$

$$\mathbb{E}\left(\widehat{W}_k^{(n)}\right)^2 \leq \frac{C_{14,2}\lambda_1\lambda_k + C_{15}\lambda_1^2\widehat{\beta}^{1-2\varepsilon}}{\lambda_1 - \lambda_k} \cdot \beta \leq (C_{14,2} + C_{15})\widehat{\beta}.$$

*Proof.* By Lemma 14 and the fact that $h_{k,n}$ forms a martingale difference sequence, we take second moment on both sides of (C.28) and obtain

$$\mathbb{E}\left(\widehat{W}_k^{(n)}\right)^2 = \sum_{i=1}^{N^*_{\beta,0.5-\varepsilon}} (1 - \beta(\lambda_1 - \lambda_k))^{2(n-i)}\mathbb{E}h_{k,i}^2 \\ + \sum_{i=N^*_{\beta,0.5-\varepsilon}+1}^{n} (1 - \beta(\lambda_1 - \lambda_k))^{2(n-i)}\mathbb{E}h_{k,i}^2 \equiv \mathrm{I} + \mathrm{II}. \quad (\mathrm{C}.30)$$

Note for I, recall from the definition in (2.3) that for each coordinate $k = 2, \ldots, d$,

$$(1 - \beta(\lambda_1 - \lambda_k))^{N^*_s} \leq \widehat{\beta}^s.$$

Hence for $n \geq N^*_{\beta,1-2\varepsilon}$

$$(1 - \beta(\lambda_1 - \lambda_k))^{2(n-N^*_{\beta,0.5-\varepsilon})} \leq \widehat{\beta}^{1-2\varepsilon},$$

and we factor this out. Applying (C.25) we have

$$\mathrm{I} \leq \widehat{\beta}^{1-2\varepsilon} \sum_{i=1}^{N^*_{\beta,0.5-\varepsilon}} (1 - \beta(\lambda_1 - \lambda_k))^{2(N^*_{\beta,0.5-\varepsilon}-i)} \mathbb{E}h_{k,i}^2 \\ \leq \widehat{\beta}^{1-2\varepsilon} \sum_{n=0}^{\infty} (1 - \beta(\lambda_1 - \lambda_k))^{2n} \cdot C_{14,1}\lambda_1^2\beta^2 \leq \frac{C_{14,1}\lambda_1^2\widehat{\beta}^{1-2\varepsilon}}{\lambda_1 - \lambda_k} \cdot \beta. \quad (\mathrm{C}.31)$$

For II in (C.30), using Lemma 14 and $(\lambda_1 - \lambda_k)\beta \leq \widehat{\beta} \leq 1$

$$\mathrm{II} \leq C_{14,2}\beta^2\left(\lambda_1\lambda_k + 4\widehat{\beta}^{1-2\varepsilon}\lambda_1^2\right) \sum_{i=N^*_{\beta,0.5-\varepsilon}+1}^{n} (1 - \beta(\lambda_1 - \lambda_k))^{2(n-i)} \\ \leq C_{14,2}\beta^2\left(\lambda_1\lambda_k + 4\widehat{\beta}^{1-2\varepsilon}\lambda_1^2\right) \sum_{n=0}^{\infty} (1 - \beta(\lambda_1 - \lambda_k))^{2n} \leq \frac{C_{14,2}\left(\lambda_1\lambda_k + 4\widehat{\beta}^{1-2\varepsilon}\lambda_1^2\right)}{\lambda_1 - \lambda_k} \cdot \beta. \quad (\mathrm{C}.32)$$

Combining (C.30), (C.31), (C.32)

$$\mathbb{E}\left(\widehat{W}_k^{(n)}\right)^2 = \mathrm{I} + \mathrm{II} \leq \frac{C_{14,2}\lambda_1\lambda_k + (C_{14,1} + 4C_{14,2}\lambda_1^2)\widehat{\beta}^{1-2\varepsilon}}{\lambda_1 - \lambda_k} \cdot \beta.$$

Setting $C_{15} = C_{14,1} + 4C_{14,2}$ in the above line proves the first inequality. The second inequality is due to $\lambda_k \leq \lambda_2$ and hence

$$\frac{\lambda_1^2\beta}{\lambda_1 - \lambda_k} \leq \frac{\lambda_1^2\beta}{\lambda_1 - \lambda_2} = \widehat{\beta},$$

completes the proof of our lemma.

$\square$

To analyze the variances in Proposition 2 we need the following probabilistic lemma, which will see its use in error estimates.

**Lemma 16.** For any random variables $X, Y$ we have
$$\mathbb{E}(X+Y)^2 \leq \mathbb{E}X^2 + 3\left(\max(\mathbb{E}X^2, \mathbb{E}Y^2) \cdot \mathbb{E}Y^2\right)^{1/2}.$$

*Proof.* We have by Cauchy-Schwarz inequality
$$\mathbb{E}(X+Y)^2 = \mathbb{E}X^2 + \mathbb{E}Y^2 + 2\mathbb{E}XY \leq \mathbb{E}X^2 + \mathbb{E}Y^2 + 2\left(\mathbb{E}X^2 \mathbb{E}Y^2\right)^{1/2}$$
$$\leq \mathbb{E}X^2 + 3\left(\max(\mathbb{E}X^2, \mathbb{E}Y^2) \cdot \mathbb{E}Y^2\right)^{1/2}.$$
□

*Proof of Proposition 2.* We set
$$\mathcal{H}_0 = \left(\mathcal{N}_v > N^*_{\beta,t}, \mathcal{N}_M > N^*_{\beta,t}\right). \tag{C.33}$$

Since (3.12) in Lemma 1 indicates
$$\mathbb{P}\left(\mathcal{N}_M \leq N^*_{\beta,t}\right) \leq 2(d+1)N^*_{\beta,t}\exp\left(-\widehat{\beta}^{-2\varepsilon}\right),$$

by letting $C_{2,P} = C_{13} \wedge 1$ we have from (C.21)
$$\mathbb{P}\left(\mathcal{H}^c_0\right) = \mathbb{P}\left(\mathcal{N}_v \leq N^*_{\beta,t}, \mathcal{N}_M > N^*_{\beta,t}\right) + \mathbb{P}\left(\mathcal{N}_M \leq N^*_{\beta,t}\right)$$
$$\leq 2(d-1)N^*_{\beta,t}\exp\left(-C_{13}\widehat{\beta}^{-2\varepsilon}\right) + 2(d+1)N^*_{\beta,t}\exp\left(-\widehat{\beta}^{-2\varepsilon}\right)$$
$$\leq 4dN^*_{\beta,t}\exp\left(-C_{2,P}\widehat{\beta}^{-2\varepsilon}\right).$$

This establishes (3.13), which leaves us to prove (3.14).

From the definition of $\mathcal{H}_0$ in (C.33) we have for each $n \in [N^*_{\beta,1-2\varepsilon}+1, N^*_{\beta,t}]$ that $h_{k,n} = f_{k,n}$. Combining this with (C.5) and (C.29), we have for $n \in [N^*_{\beta,0.5-\varepsilon}+1, N^*_{\beta,t}]$
$$\left(U^{(n)}_k\right)^2 \mathbf{1}_{\mathcal{H}_0} = \left(W^{(n)}_k\right)^2.$$

Taking expectation
$$\mathbb{E}\left[\left(U^{(n)}_k\right)^2; \mathcal{H}_0\right] = \mathbb{E}\left(W^{(n)}_k\right)^2. \tag{C.34}$$

To compute the right hand, we have from (C.29)
$$W^{(n)}_k = (1-\beta(\lambda_1-\lambda_k))^n U^{(0)}_k + \widehat{W}^{(n)}_k + \sum_{i=1}^n (1-\beta(\lambda_1-\lambda_k))^{n-i} E_{k,i-1}.$$

Note from Lemma 9(i) we have $|E_{k,n-1}| \leq C_9 M^2 \beta^2 = C_9 \widehat{\beta}^{-4\varepsilon}\beta^2$, so
$$\left|\sum_{i=1}^n (1-\beta(\lambda_1-\lambda_k))^{n-i} E_{k,i-1}\right| \leq \sum_{i=1}^n (1-\beta(\lambda_1-\lambda_k))^{n-i} |E_{k,i-1}|$$
$$\leq C_9 \widehat{\beta}^{-4\varepsilon}\beta^2 \sum_{i=1}^n (1-\beta(\lambda_1-\lambda_k))^{n-i} \leq \frac{C_9 \lambda^2_1 \widehat{\beta}^{-4\varepsilon}\beta}{\lambda_1-\lambda_k}.$$

Combine the last three displays and Lemma 16, using (C.2) and the fact that $\widehat{\beta} \leq \widehat{\beta}^{1-2\varepsilon} \leq \widehat{\beta}^{0.5-4\varepsilon} \leq 1$ we conclude from Lemma 16 and Lemma 15

$$
\begin{aligned}
\mathbb{E}\left(W_k^{(n)}\right)^2 &= (1-\beta(\lambda_1-\lambda_k))^{2n}(U_k^{(0)})^2 + \mathbb{E}\left(\widehat{W}_k^{(n)} + \sum_{i=1}^n (1-\beta(\lambda_1-\lambda_k))^{n-i} E_{k,i-1}\right)^2 \\
&\leq (1-\beta(\lambda_1-\lambda_k))^{2n}\left(U_k^{(0)}\right)^2 \\
&\quad + \frac{C_{14,2}\lambda_1\lambda_k + C_{15}\lambda_1^2\widehat{\beta}^{1-2\varepsilon}}{\lambda_1-\lambda_k}\cdot \beta + 3\left(C_{14,2}+C_{15}\right)^{0.5}\widehat{\beta}^{0.5}\cdot \frac{C_9\lambda_1^2\widehat{\beta}^{-4\varepsilon}}{\lambda_1-\lambda_k}\beta \\
&\leq (1-\beta(\lambda_1-\lambda_k))^{2n}\left(U_k^{(0)}\right)^2 \\
&\quad + \frac{C_{14,2}\lambda_1\lambda_k + \left[C_{15}+3(C_{14,2}+C_{15})^{0.5}C_9\right]\lambda_1^2\cdot \widehat{\beta}^{0.5-4\varepsilon}}{\lambda_1-\lambda_k}\cdot \beta.
\end{aligned}
$$

Setting $C'_{2,P} = \max\left(C_{14,2}, C_{15}+3(C_{14,2}+C_{15})^{0.5}C_9\right)$ in the last line above, and noting the relation in (C.34) completes proof of (3.14), and hence the proposition.

□

### C.3 Proof of Proposition 3

Let the univariate function $f(x) = x(1-x^2)^{-1/2}$. First we analyze the increment of

$$f(v_1^{(n)}) = \frac{v_1^{(n)}}{\left[1-(v_1^{(n)})^2\right]^{1/2}}.$$

We conclude an incremental type of lemma as follows.

**Lemma 17.** Assume all conditions in Theorem 1. Then there exists a random variable $R_n$ satisfying that on the event $(\mathcal{N}_M > n, \mathcal{N}_c \geq n)$ we have

$$|R_n| \leq C_{17}M^2\beta^2, \quad a.s., \tag{C.35}$$

where $\mathcal{N}_M$ was defined earlier in (3.6), such that the increment of $U_k^{(n)}$ has

$$f(v_1^{(n)}) - f(v_1^{(n-1)}) = \frac{\beta}{(1-v_1^2)^{3/2}}\left[(\mathbf{v}^\top \boldsymbol{Y})Y_1 - v_1(\mathbf{v}^\top \boldsymbol{Y})^2\left(1+\frac{\beta}{2}\|\boldsymbol{Y}\|^2\right)\right] + R_n. \tag{C.36}$$

*Proof.* The derivative of $f(x)$ is

$$f'(x) = \frac{\mathrm{d}}{\mathrm{d}x}\left[\frac{x}{(1-x^2)^{1/2}}\right] = \frac{1}{(1-x^2)^{3/2}}$$

Taylor's theorem with the Lagrange form of the remainder gives

$$f(\widehat{\mathbf{v}}) - f(\mathbf{v}) = f'(v_1)(\widehat{v}_1 - v_1) + R'_n = \frac{\widehat{v}_1 - v_1}{\left(1-v_1^2\right)^{3/2}} + R'_n, \tag{C.37}$$

where there is a constant $C \leq 12$ that bound the second derivative $f''(x)$ for $x \in [-1/\sqrt{2}, 1/\sqrt{2}]$, and hence
$$|R'_n| \leq \frac{C}{2}(\widehat{v}_1 - v_1)^2 \leq 6C_{4,2}^2 M^2 \beta^2,$$
where we applied (A.2). Note (A.1) in Proposition 4 implies
$$\widehat{v}_1 - v_1 = \beta \left[ (\mathbf{v}^\top \boldsymbol{Y}) Y_k - v_k (\mathbf{v}^\top \boldsymbol{Y})^2 \left(1 + \frac{\beta}{2} \|\boldsymbol{Y}\|^2 \right)\right] + Q_{1,n}. \tag{C.38}$$

Combining (C.37) and (C.38) we obtain
$$f(\widehat{\mathbf{v}}) - f(\mathbf{v}) = f'(v_1)(\widehat{v}_1 - v_1) + R'_n$$
$$= \frac{\beta}{(1-v_1^2)^{3/2}} \left[ (\mathbf{v}^\top \boldsymbol{Y}) Y_k - v_k (\mathbf{v}^\top \boldsymbol{Y})^2 \left(1 + \frac{\beta}{2} \|\boldsymbol{Y}\|^2 \right)\right] + \frac{Q_{1,n}}{(1-v_1^2)^{3/2}} + R'_n.$$

Setting $R_n = Q_{1,n} \left(1 - v_1^2\right)^{-3/2} + R'_n$ completes the proof of (C.36). $\square$

Let $\alpha = 1 + (\beta/2) d \cdot \lambda_1 \widehat{\beta}^{-2\varepsilon}$, and hence on $(\mathcal{N}_M > n)$ one has $1 + (\beta/2)\|\boldsymbol{Y}\|^2 \leq \alpha$. We define the coupled process $U^{(n)}$ as $U^{(0)} = f(v_1^{(0)})$ and recursively
$$U^{(n)} = U^{(n-1)} + \frac{\beta}{(1-v_1^2)^{3/2}} \left[ (\mathbf{v}^\top \boldsymbol{Y}) Y_1 - \alpha v_1 (\mathbf{v}^\top \boldsymbol{Y})^2 \right] + R_n.$$

When $v_1^{(n-1)} > 0$ the increment of $f(\mathbf{v}^{(n)})$ is lower-bounded by that of $U^{(n)}$, i.e.,
$$f(\mathbf{v}^{(n)}) - f(\mathbf{v}^{(n-1)}) \geq U^{(n)} - U^{(n-1)}. \tag{C.39}$$

We turn to analyze the increment $U^{(n)} - U^{(n-1)}$ by first calculating the expectation

**Lemma 18.** On $(\mathcal{N}_c \geq n) \in \mathcal{F}_{n-1}$ we have
$$\mathbb{E}\left[ \frac{\beta}{(1-v_1^2)^{3/2}} \left( (\mathbf{v}^\top \boldsymbol{Y}) Y_1 - \alpha v_1 (\mathbf{v}^\top \boldsymbol{Y})^2 \right) \,\Big|\, \mathcal{F}_{n-1} \right] = \beta \cdot \frac{\lambda_1 - \alpha \mathbf{v}^\top \boldsymbol{\Lambda} \mathbf{v}}{1 - v_1^2} \cdot \frac{v_1}{(1-v_1^2)^{1/2}}. \tag{C.40}$$

*Proof.* Note for all $\mathbf{v} \in \mathcal{S}^{d-1}$, we have $\mathbb{E}[(\mathbf{v}^\top \boldsymbol{Y}) \mid \mathcal{F}_{n-1}] = \lambda_1 v_1$ and
$$-\mathbb{E}\left[\alpha v_1 (\mathbf{v}^\top \boldsymbol{Y})^2 \mid \mathcal{F}_{n-1}\right] = -\alpha v_1 \mathbf{v}^\top \mathbb{E}\left[\boldsymbol{Y} \boldsymbol{Y}^\top \mid \mathcal{F}_{n-1}\right] \mathbf{v} = -\alpha v_1 \mathbf{v}^\top \boldsymbol{\Lambda} \mathbf{v}.$$

Adding up the above two equations proves (C.40). $\square$

Now from (C.36)
$$U^{(n)} - U^{(n-1)} = \frac{\beta}{(1-v_1^2)^{3/2}} \left[ (\mathbf{v}^\top \boldsymbol{Y}) Y_1 - \alpha v_1 (\mathbf{v}^\top \boldsymbol{Y})^2 \right] + R_n.$$

Let
$$e_n = \frac{\beta}{(1-v_1^2)^{3/2}} \left[ (\mathbf{v}^\top \boldsymbol{Y}) Y_1 - \alpha v_1 (\mathbf{v}^\top \boldsymbol{Y})^2 \right] - \beta \cdot \frac{\lambda_1 - \alpha \mathbf{v}^\top \boldsymbol{\Lambda} \mathbf{v}}{1 - v_1^2} \cdot \frac{v_1}{(1-v_1^2)^{1/2}},$$

then $\bar{e}_n$ forms a martingale difference sequence, and from (C.40) in Lemma 18

$$U^{(n)} = e_n + \left(1 + \beta \cdot \frac{\lambda_1 - \alpha \mathbf{v}^\top \mathbf{\Lambda} \mathbf{v}}{1 - v_1^2}\right) \cdot U^{(n-1)} + R_n, \tag{C.41}$$

where since $d\widehat{\beta}^{1-2\varepsilon} \leq 1$

$$\frac{\lambda_1 - \alpha \mathbf{v}^\top \mathbf{\Lambda} \mathbf{v}}{1 - v_1^2} \geq \left(\alpha - \frac{1}{2}\right)(\lambda_1 - \lambda_2) \geq \frac{\lambda_1 - \lambda_2}{2}. \tag{C.42}$$

Furthermore let

$$\bar{e}_n = e_n \mathbf{1}_{\mathcal{N}_c \geq n},$$

and

$$\bar{R}_n = R_n \mathbf{1}_{\mathcal{N}_M > n, \mathcal{N}_c \geq n},$$

let $\bar{U}^{(0)} = 0$, and we define a coupled process $\bar{U}^{(n)}$ as

$$\bar{U}^{(n)} = \bar{e}_n + \left(1 + \beta \cdot \frac{\lambda_1 - \alpha \mathbf{v}^\top \mathbf{\Lambda} \mathbf{v}}{1 - v_1^2}\right) \cdot \bar{U}^{(n-1)} + \bar{R}_n, \tag{C.43}$$

we conclude the following coupling lemma that is in the same fashion as Lemma 11.

**Lemma 19.** For each $n \geq 0$ we have $\bar{U}^{(n)} = U^{(n)}$ on the event $(\mathcal{N}_M > n, \mathcal{N}_c \geq n)$.

*Proof.* By definition the equality is valid for $n = 0$. Suppose the lemma is valid for $n - 1$, then on the event $(\mathcal{N}_M > n, \mathcal{N}_c \geq n)$ which is a subset of $(\mathcal{N}_M > n - 1, \mathcal{N}_c \geq n - 1)$ we have

$$\bar{U}^{(n-1)} = U^{(n-1)},$$

and also $e_n = \bar{e}_n$ and $R_n = \bar{R}_n$. From this we immediately obtain from (C.41) and (C.43) that $\bar{U}^{(n)} = U^{(n)}$. □

Let for each $n$

$$P_{-n} = \prod_{i=0}^{n-1} \left(1 + \beta \cdot \frac{\lambda_1 - \alpha \mathbf{v}^{(i)\top} \mathbf{\Lambda} \mathbf{v}^{(i)}}{1 - (v_1^{(i)})^2}\right)^{-1}. \tag{C.44}$$

Then $P_{-n}$ is a $\mathcal{F}_{n-1}$-measurable random variable, and from (C.42)

$$P_{-n} \leq \left(1 + \beta \cdot \frac{\lambda_1 - \lambda_2}{2}\right)^{-n}.$$

Let

$$E_{n-1} = \mathbb{E}\left[\bar{R}_n \mid \mathcal{F}_{n-1}\right],$$

and let

$$f_n = \bar{e}_n + \bar{R}_n - E_{n-1},$$

then we have the following expression

**Lemma 20.** The sequence generated by (C.43) has representation

$$P_{-n}\bar{U}^{(n)} = \bar{U}^{(0)} + \sum_{i=1}^{n} P_{-i}\left(f_i + E_{i-1}\right), \qquad (C.45)$$

*Proof.* (C.43) is just

$$\bar{U}^{(n)} = f_n + \left(1 + \beta \cdot \frac{\lambda_1 - \alpha \mathbf{v}^{(n-1)\,\top} \mathbf{\Lambda} \mathbf{v}^{(n-1)}}{1 - v_1^2}\right) \cdot \bar{U}^{(n-1)} + E_{n-1}.$$

Multiplying both sides by $P_{-n}$ we have

$$P_{-n}\bar{U}^{(n)} = P_{-n}f_n + P_{-n}\left(1 + \beta \cdot \frac{\lambda_1 - \alpha \mathbf{v}^{(n-1)\,\top} \mathbf{\Lambda} \mathbf{v}^{(n-1)}}{1 - v_1^2}\right) \cdot \bar{U}^{(n-1)} + P_{-n}E_{n-1}$$

$$= P_{-n}f_n + P_{-(n-1)}\bar{U}^{(n-1)} + P_{-n}E_{n-1}.$$

Iteratively applying this and noting $P_0 = 1$

$$P_{-n}\bar{U}^{(n)} = \bar{U}^{(0)} + \sum_{i=1}^{n} P_{-i}(f_i + E_{i-1}),$$

concluding (C.45) and hence the lemma. $\square$

**Lemma 21.** For all $y$ such that

$$y \leq \frac{C'_*}{2} \cdot \widehat{\beta}^{-0.5+\varepsilon}, \qquad (C.46)$$

We have for each fixed $n$

$$\mathbb{P}\left(\left|\sum_{i=1}^{n} P_{-i}f_i\right| \geq y\widehat{\beta}^{0.5-\varepsilon}\right) \leq 2\exp\left(-C_{21}y^2\widehat{\beta}^{-2\varepsilon}\right). \qquad (C.47)$$

*Proof.* (i) The proof is analogous to the one of Lemma 12. We first establish that $f_n$ forms a martingale difference sequence and is subexponential. To show this, first note

$$f_n = \bar{e}_n + \bar{R}_n - \mathbb{E}[\bar{R}_n \mid \mathcal{F}_{n-1}].$$

Note from (C.35) we have $|\bar{R}_n| \leq C_{17}M^2\beta^2, a.s.$ Therefore

$$\left\|\bar{R}_n - \mathbb{E}[\bar{R}_n \mid \mathcal{F}_{n-1}]\right\|_{\psi_1} \leq 2C_{17}M^2\beta^2.$$

Hence from the definition of $\mathcal{S}_1$ and $\mathcal{N}_c$

$$\|\bar{e}_n\|_{\psi_1} \leq \left\|\mathbf{1}_{\mathcal{N}_c \geq n} \frac{2\beta}{(1-v_1^2)^{3/2}}\left[(\mathbf{v}^\top \mathbf{Y})Y_1 - \alpha v_1 (\mathbf{v}^\top \mathbf{Y})^2\right]\right\|_{\psi_1}$$

$$\leq 4\sqrt{2}\beta\left(\|(\mathbf{v}^\top \mathbf{Y})Y_1\|_{\psi_1} + \|\alpha v_1 (\mathbf{v}^\top \mathbf{Y})^2\|_{\psi_1}\right).$$

However from Assumption 1, Lemma 6 implies
$$\|(\mathbf{v}^\top \boldsymbol{Y})Y_1\|_{\psi_1} \le \|\mathbf{v}^\top \boldsymbol{Y}\|_{\psi_2}^2 + \|Y_1\|_{\psi_2}^2 \le 2\lambda_1,$$
and
$$\|\alpha v_1(\mathbf{v}^\top \boldsymbol{Y})^2\|_{\psi_1} \le \|(\mathbf{v}^\top \boldsymbol{Y})^2\|_{\psi_1} \le \|\mathbf{v}^\top \boldsymbol{Y}\|_{\psi_2}^2 \le 2\lambda_1.$$

We have
$$\|\bar{e}_n\|_{\psi_1} \le 16\sqrt{2}\lambda_1\beta.$$

Therefore from $\lambda_1\beta \le \widehat{\beta} \le \widehat{\beta}^{4\varepsilon} \implies M^2\beta^2 = \lambda_1^2\widehat{\beta}^{-4\varepsilon}\beta^2 \le \lambda_1\beta$

$$\|f_n\|_{\psi_1} \le 16\sqrt{2}\lambda_1\beta + 2C_{17}M^2\beta^2 \le (16\sqrt{2} + 2C_{17})\lambda_1\beta.$$

Therefore if $K$ is chosen as $C_{7,1}C'_*\lambda_1\beta$ where $C'_* = 16\sqrt{2} + 2C_{17}$, applying Lemma 7 we have for $|t| \le K^{-1}$
$$\mathbb{E}[\exp(tf_n) \mid \mathcal{F}_{n-1}] \le \exp\left(\frac{t^2 K^2}{2}\right). \tag{C.48}$$

(ii) We now come to prove a concentration inequality for $\sum_{i=1}^n P_{-i}f_i$. For $i = 1, 2, \ldots, n$
$$P_{-i} \le \left(1 + \beta \cdot \frac{\lambda_1 - \lambda_2}{2}\right)^{-i} \le 1.$$

Set $a_i = P_{-i}$ as the $i$th coordinate of $\mathbf{a}$ for $i = 1, \ldots, n$, and $f_i$ in the place of $\xi_i$, then
$$\|\mathbf{a}\|_2^2 = \sum_{i=1}^n P_{-i}^2 \le \sum_{i=1}^n \left(1 + \beta \cdot \frac{\lambda_1 - \lambda_2}{2}\right)^{-2i} \le \frac{1}{\beta(\lambda_1 - \lambda_2)}, \tag{C.49}$$

Applying Lemma 8 with
$$\|\mathbf{a}\|_\infty \le 1, \qquad \|\mathbf{a}\|_2 \le (\lambda_1 - \lambda_2)^{-0.5}(2\beta)^{-0.5}, \qquad K = C'_*\lambda_1\beta, \qquad z = y\widehat{\beta}^{0.5-\varepsilon},$$

from the definition of $\widehat{\beta}$ we have
$$\frac{z^2}{K^2\|\mathbf{a}\|_2^2} \ge \frac{y^2\widehat{\beta}^{1-2\varepsilon}}{C'^2_*\lambda_1^2\beta^2 \cdot (\lambda_1 - \lambda_2)^{-1}\beta^{-1}} \tag{C.50}$$
$$= C'^{-2}_* y^2\widehat{\beta}^{1-2\varepsilon} \cdot (\lambda_1 - \lambda_2)\lambda_1^{-2} \cdot \beta^{-1} = C'^{-2}_* y^2\widehat{\beta}^{-2\varepsilon},$$

and from $\lambda_1\beta \le \widehat{\beta}$
$$\frac{z}{K\|\mathbf{a}\|_\infty} \ge \frac{y\widehat{\beta}^{0.5-\varepsilon}}{C'_*\lambda_1\beta \cdot 1} = C'^{-1}_* y \cdot \lambda_1^{-1} \cdot \widehat{\beta}^{0.5-\varepsilon}\beta^{-1} \ge C'^{-1}_* y \cdot \widehat{\beta}^{-0.5-\varepsilon}. \tag{C.51}$$

Whenever (C.46) holds we have $C'^{-2}_* y^2\widehat{\beta}^{-2\varepsilon} \le C'^{-1}_* y \cdot \widehat{\beta}^{-0.5-\varepsilon}$, the minimum of (C.50) and (C.51) is
$$\min\left(\frac{z^2}{K^2\|\mathbf{a}\|_2^2}, \frac{z}{K\|\mathbf{a}\|_\infty}\right) \ge \min\left(C'^{-2}_* y^2\widehat{\beta}^{-2\varepsilon}, C'^{-1}_* y \cdot \widehat{\beta}^{-0.5-\varepsilon}\right) = C'^{-2}_* y^2\widehat{\beta}^{-2\varepsilon}.$$

Hence using (C.49) and the concentration inequality as in Lemma 8 we proves

$$\mathbb{P}\left(\left|\sum_{i=1}^{n} P_{-i} f_i\right| \geq y \widehat{\beta}^{0.5-\varepsilon}\right) \leq 2 \exp\left(-C_8 \cdot C_*'^{-2} \cdot y^2 \widehat{\beta}^{-2\varepsilon}\right). \tag{C.52}$$

Letting $C_{21} = C_8 \cdot C_*'^{-2}$ in (C.52) proves the lemma. $\square$

With the concentration inequality Lemma 21 at hand, we are ready for the proof of Proposition 3.

*Proof of Proposition 3.* (i) Define

$$\mathcal{E} = \left(\sup_{n \leq N_\beta^o(c^*)} \left|P_{-n} \bar{U}^{(n)} - U^{(0)}\right| \leq 2C_{17} \widehat{\beta}^{0.5-\varepsilon}\right) \bigcap \left(\mathcal{N}_M > N_\beta^o(c^*)\right). \tag{C.53}$$

First we want to estimate $\mathbb{P}(\mathcal{E}^c)$. Note from (C.45) in Lemma 20 that

$$P_{-n} \bar{U}^{(n)} - U^{(0)} = \sum_{i=1}^{n} P_{-i}(f_i + E_{i-1}),$$

and since $\varepsilon \leq 1/6$, $\widehat{\beta}^{1-4\varepsilon} \leq \widehat{\beta}^{0.5-\varepsilon}$ and we have

$$\sum_{i=1}^{n} P_{-i}|E_{i-1}| \leq \sum_{i=1}^{n} (1 + \beta(\lambda_1 - \lambda_2))^{-i} \cdot C_{17} M^2 \beta^2$$

$$\leq \frac{C_{17} M^2 \beta}{\lambda_1 - \lambda_2} = C_{17} \widehat{\beta}^{1-4\varepsilon} \leq C_{17} \widehat{\beta}^{0.5-\varepsilon}.$$

Taking $y = C_{17}$ in (C.47) we have for each $n \geq 1$

$$\mathbb{P}\left(|P_{-n} \bar{U}^{(n)} - U^{(0)}| \geq 2C_{17} \widehat{\beta}^{0.5-\varepsilon}\right) \leq \mathbb{P}\left(\left|\sum_{i=1}^{n} P_{-i} f_i\right| \geq C_{17} \widehat{\beta}^{0.5-\varepsilon}\right)$$

$$\leq 2 \exp\left(-C_{21} C_{17}^2 \widehat{\beta}^{-2\varepsilon}\right),$$

and hence from (3.12) we have

$$\mathbb{P}(\mathcal{E}^c) \leq \sum_{n=1}^{N_\beta^o(c^*)} \mathbb{P}\left(\left|P_{-n} \bar{U}^{(n)} - U^{(0)}\right| > 2C_{17} \widehat{\beta}^{0.5-\varepsilon}\right) + \mathbb{P}\left(\mathcal{N}_M \leq N_\beta^o(c^*)\right) \tag{C.54}$$

$$\leq 2N_\beta^o(c^*) \exp\left(-C_{21} C_{17}^2 \widehat{\beta}^{-2\varepsilon}\right) + 2(d+1) N_\beta^o(c^*) \exp\left(-\widehat{\beta}^{-2\varepsilon}\right).$$

(ii) On the event $\mathcal{E} \cap \left(\mathcal{N}_c > N_\beta^o(c^*)\right)$ we have from (C.53) and the coupling Lemma 19 that

$$\sup_{n \leq N_\beta^o(c^*)} \left|P_{-n} U^{(n)} - U^{(0)}\right| \leq 2C_{17} \widehat{\beta}^{0.5-\varepsilon}.$$

Also, the initial condition in Theorem 1 is equivalent to

$$|U^{(0)}| = |\cot \angle(\mathbf{u}^{(0)}, \mathbf{u}^*)| \geq (c^*d)^{-0.5}.$$

Therefore from (2.4) we know by choosing $b_5 = (16C_{17})^{-1}$ in (3.16), i.e. $2C_{17}\widehat{\beta}^{0.5-\varepsilon} \leq (1/2)(c^*d)^{-0.5}$

$$\left|U^{(N_\beta^o(c^*))}\right| \geq \left(P_{-N_\beta^o(c^*)}\right)^{-1}\left(|U^{(0)}| - 2C_{17}\widehat{\beta}^{0.5-\varepsilon}\right)$$
$$\geq \left(1 + \beta \cdot \frac{\lambda_1 - \lambda_2}{2}\right)^{N_\beta^o(c^*)}\left((c^*d)^{-0.5} - 2C_{17}\widehat{\beta}^{0.5-\varepsilon}\right)$$
$$\geq \left(1 + \beta \cdot \frac{\lambda_1 - \lambda_2}{2}\right)^{N_\beta^o(c^*)} \cdot \frac{1}{2}(c^*d)^{-0.5} \geq 1.$$

By taking logarithm in the last display and noticing $0 \leq t \leq 1/2$ both the following hold:

$$\frac{\log(1 + t/2)}{-\log(1 - t)} \geq \frac{1 - t}{2}, \qquad \log(1 - t) \geq -2t,$$

we have for $d\widehat{\beta}^{1-2\varepsilon} \leq 1/c^* \leq 1$

$$\log|U^{(N_\beta^o(c^*))}| \geq N_\beta^o(c^*)\log(1 + \beta(\lambda_1 - \lambda_2)/2) + \log\left((c^*d)^{-0.5} - 2C_{17}\widehat{\beta}^{0.5-\varepsilon}\right)$$
$$\geq \log(4c^*d)\frac{\log(1 + \beta(\lambda_1 - \lambda_2)/2)}{-\log(1 - \beta(\lambda_1 - \lambda_2))}$$
$$\quad - 0.5\log(c^*d) + \log\left(1 - 2C_{17}(c^*d\widehat{\beta}^{1-2\varepsilon})^{0.5}\right)$$
$$\geq 0.5\log(4c^*d)(1 - \beta(\lambda_1 - \lambda_2)) - 0.5\log(c^*d) - 4C_{17}(c^*d\widehat{\beta}^{1-2\varepsilon})^{0.5}$$
$$\geq 0.5\ln 4 - 0.5\log(4c^*d)\widehat{\beta} - 4C_{17}(c^*d\widehat{\beta}^{1-2\varepsilon})^{0.5}.$$

It is straightforward to verify that the last line above is nonnegative whenever $c^*d\widehat{\beta}^{1-2\varepsilon} \leq \ln^2 2/(4C_{17} + 2)^2$, which gives (3.16), and hence by exponentiation we obtain $|U^{(N_\beta^o(c^*))}| \geq 1$, and hence from the increment relation C.39 we have $|f(vb^{(n)})| \geq |U^{(N_\beta^o(c^*))}| \geq 1$

From the definition of $\mathcal{N}_c$ in (3.15) the above proves

$$\mathcal{E} \cap \left(\mathcal{N}_c > N_\beta^o(c^*)\right) \subseteq \left(\mathcal{N}_c \leq N_\beta^o(c^*)\right) \iff \mathcal{E} \subseteq \left(\mathcal{N}_c \leq N_\beta^o(c^*)\right).$$

Combining the last display with (C.54) and letting $C_{3,P} = (C_{21}C_{17}^2) \wedge 1$ indicate

$$\mathbb{P}\left(\mathcal{N}_c \leq N_\beta^o(c^*)\right) \geq \mathbb{P}(\mathcal{E}) \geq 1 - 2(d+2)N_\beta^o(c^*)\exp\left(-C_{3,P}\widehat{\beta}^{-2\varepsilon}\right),$$

proving (3.17).

$\square$

# D  Proofs of Auxillary Lemmas

## D.1  Proof of Lemma 1

*Proof of Lemma 1.* We have from union bound

$$\mathbb{P}\left(\mathcal{N}_M \leq N\right) \leq \sum_{n=1}^N \mathbb{P}\left(\max_{1\leq k \leq d}\left|Y_k^{(n)}\right|, \left|\mathbf{v}^{(n-1)\top}\mathbf{Y}^{(n)}\right| \geq M^{1/2}\right)$$

$$\leq \sum_{n=1}^N \left(\sum_{k=1}^d \mathbb{P}\left(\left|Y_k^{(n)}\right| \geq M^{1/2}\right) + \mathbb{P}\left(\left|\mathbf{v}^{(n-1)\top}\mathbf{Y}^{(n)}\right| \geq M^{1/2}\right)\right)$$

Since $Y_k^{(n)}$, $k=1,2,\ldots,d$ and $\mathbf{v}^{(n-1)\top}\mathbf{Y}^{(n)}$ has subgaussian norm $\leq \lambda_1^{1/2}$, Lemma 4(i) implies for an arbitrary subgaussian random variable $Z$

$$\mathbb{P}(|Z| \geq M^{1/2}) \leq 2\exp\left(-M/\|Z\|_{\psi_2}^2\right),$$

we have by choosing $M = \lambda_1 \widehat{\beta}^{-2\varepsilon}$ in (3.11) that the probability is no greater than

$$2(d+1)N\exp\left(-M/\lambda_1\right) = 2(d+1)N\exp\left(-\widehat{\beta}^{-2\varepsilon}\right).$$

□

## D.2  Proof of Lemma 2

To prove Lemma 2 we need to use a classical result of multivariate statistics.

**Lemma 22.** *If $\mathbf{X} = (X_1,\ldots,X_d)^\top$ is uniformly generated from $\mathcal{S}^{d-1}$ in $\mathbb{R}^d$, the marginal density of $X_1$ is*

$$\frac{\mathrm{d}}{\mathrm{d}x}\mathbb{P}(X_1 \leq x) = \frac{\omega_{d-1}}{\omega_d}(1-x)^{(d-3)/2} \cdot 1_{[-1,1]}(x),$$

*where $\omega_{d_*} \equiv 2\pi^{d_*/2}[\Gamma(d_*/2)]^{-1}$ is the surface area of the unit sphere $\mathcal{S}^{d_*-1}$ in $\mathbb{R}^{d_*}$.*

Proof of Lemma 22 can be found in (Muirhead, 2005, pp. 147) and is hence omitted.

*Proof of Lemma 2.* Let $\theta_0 = \angle(\mathbf{u}^{(0)}, \mathbf{u}^*)$. Since $\cos\theta_0$ is distributed the same as $X_1$ in Lemma 22 we use the density formula there to conclude

$$\mathbb{P}\left(\tan^2\theta_0 > C^*\delta^{-2}d\right) = \mathbb{P}\left(|\cos\theta_0| < \left(1+C^*\delta^{-2}d\right)^{1/2}\right)$$
$$\leq 4\frac{\omega_{d-1}}{\omega_d}(C^*\delta^{-2}d)^{-1/2} \leq \delta,$$

where $C^* = \sup_{d\geq 1} 16(\omega_{d-1}/\omega_d)^2 d^{-1}$ is finite due to the property of gamma function. The numerical value of $C^*$ is approximately 2.56.

□